\input ./style/arxiv-general.cfg
\documentclass[aap,MSNbibl,secthm,seceqn,dvips]{arximspdf}
\makeatletter
   \@ifpackageloaded{graphicx}{}{\usepackage{graphicx}}
\makeatother
\usepackage{mathbh}

\doi{10.1214/15-AAP1126}% Updated by VTEXPTS2LaTeX.exe, 31.07.2015
\volume{26}
\issue{3}
\pubyear{2016}
\firstpage{1581}
\lastpage{1619}
\docsubty{FLA}

\makeatletter
\newcommand{\rrvert}{\vert}
\newcommand{\llvert}{\vert}
\newcommand{\eqref}[1]{(\ref{#1})}
\newtheorem{Theorem}{Theorem}[section]
\newtheorem{Proposition}[Theorem]{Proposition}
\newtheorem{Lemma}[Theorem]{Lemma}
\newtheorem{Corollary}[Theorem]{Corollary}
\newproclaim{Definition}[Theorem]{Definition}
\newproclaim{Assumption}[Theorem]{Assumption}
\newproclaim{remark}[Theorem]{Remark}
\newcommand{\BBone}{\mathbh{1}}
\makeatother

\begin{document}
\begin{frontmatter}

\title{Estimation for stochastic damping Hamiltonian systems under
partial observation. III.~Diffusion~term}
\runtitle{Estimation for kinetic equations}

\begin{aug}
% Corresponding author: Cl\'{e}mentine Prieur - cprieur@imag.fr% Updated by VTEXPTS2LaTeX.exe, 06.08.2015 07:48
%Updated by VTEXPTS2LaTeX.exe, 05.08.2015 15:35
%Updated by VTEXPTS2LaTeX.exe, 31.07.2015 09:15
\author[A]{\fnms{Patrick}~\snm{Cattiaux}\thanksref{m1}\ead[label=e1]{cattiaux@math.univ-toulouse.fr}},
\author[B]{\fnms{Jos\'e R.}~\snm{Le\'on}\thanksref{m2,T2,T3}\ead[label=e2]{jose.leon@ciens.ucv.ve}}
\and
\author[C]{\fnms{Cl\'ementine}~\snm{Prieur}\thanksref{m3,T2}\corref{}\ead[label=e3]{clementine.prieur@imag.fr}}
\runauthor{P. Cattiaux, J. R. Le\'on and C. Prieur}
\affiliation{Universit\'e de Toulouse\thanksmark{m1}, Universidad
Central de Venezuela\thanksmark{m2} and Universit\'e Grenoble
Alpes\thanksmark{m3}}
%\dedicated{}
\address[A]{P. Cattiaux\\
Institut de Math\'ematiques de Toulouse\\
Universit\'e de Toulouse\\
Toulouse\\
France\\
\printead{e1}}
\address[B]{J. R. Le\'on\\
Escuela de Matem\'atica, Facultad de Ciencias\\
Universidad Central de Venezuela\\
Caracas\\
Venezuela\\
\printead{e2}}
\address[C]{C. Prieur\\
Universit\'e Grenoble Alpes \\
\quad and Inria Grenoble Rh\^one-Alpes\\
Laboratoire Jean Kuntzmann, Inria/AIRSEA\\
Grenoble\\
France\\
\printead{e3}}
\end{aug}
\thankstext{T2}{Supported in part by the project ECOS NORD/FONACIT
(France-Venezuela) V12M01 \textit{Equations diff\'erentielles
stochastiques et analyse de sensibilit\'e : applications
environnementales}.}
\thankstext{T3}{Supported in part by the Inria International Chairs program.}

% HISTORY:
%
\received{\smonth{7} \syear{2014}}% Updated by VTEXPTS2LaTeX.exe,
%31.07.2015 09:15
%
\revised{\smonth{5} \syear{2015}}% Updated by VTEXPTS2LaTeX.exe,
%31.07.2015 09:15

% ABSTRACT
%
\begin{abstract}
This paper is the third part of our study started with
Cattiaux, Le\'{o}n  and Prieur
[\textit{Stochastic Process. Appl.} \textbf{124} (2014) 1236--1260;
\textit{ALEA Lat. Am. J. Probab. Math. Stat.} \textbf{11} (2014)  359--384].
For some ergodic Hamiltonian systems, we obtained a
central limit theorem for a nonparametric estimator of the invariant
density [\textit{Stochastic Process. Appl.} \textbf{124} (2014) 1236--1260]
and of the drift term [\textit{ALEA Lat. Am. J. Probab. Math. Stat.} \textbf{11} (2014)  359--384], under partial
observation (only the positions are observed). Here, we obtain
similarly a central limit theorem for a nonparametric estimator of the
diffusion term.
\end{abstract}

% KEYWORDS
\end{frontmatter}

%s1 #&#
\section{Introduction}\label{sec1}
In this article, we consider the estimation, using data sampled at high
frequency, of the local variance or diffusion term $\sigma(\cdot
,\cdot)$ in the system $ (Z_t:=(X_t,Y_t) \in\mathbb{R}^{2d},
t \geq0 )$ governed by the following It\^{o} stochastic
differential equation:
%
%e1.1 #&#
\begin{equation}
\label{WuGirsanov} \cases{ dX_t = Y_t \,dt,
\cr
\vspace*{2pt}
dY_t = \sigma(X_t,Y_t) \,dW_t-
\bigl(c(X_t,Y_t)Y_t+\nabla V(X_t)
\bigr) \,dt. }
\end{equation}

The function $c$ is called the damping force and $V$ the potential,
$\sigma$ is the diffusion term and $W$ a standard Brownian motion.

The problem of estimating the diffusion term, sometimes called
volatility, in a model of diffusion has a somewhat long history and has
a lot of motivations, in particular in the analysis of financial or
neuronal data.

The beginning of the story takes place at the end of the eighties of
the last century. The first and seminal articles were written by
\cite{DAF,DF,DOH} and \cite{GCJ}. The method
generally used is the central limit theorem for martingales. Recently,
an excellent survey introducing the subject and giving some important
recent references was written by \cite{PV}. In that work, the authors
give some insights about the methods of proof of the limit theorems and
recall also the existence of some goodness of fit tests useful in
financial studies. This article also mentioned the names of those
linked to the development in this area. They are, among others,
\cite{BS,JP} and \cite{BN}. The second of the last cited works
contains a deep study for the asymptotic behavior of discrete
approximations of stochastic integrals; it is thus in tight
relationship with the estimation of the diffusion term.

The present article is the continuation of two previous works by the
authors: \cite{CLP2} and \cite{CLP}. In the first one, we tackled the
problem of estimating the invariant density of the system (\ref{WuGirsanov}) and in the second one the estimation of the drift term
$(x,y) \mapsto b(x,y)=-c(x,y)y+\nabla V(x)$ was studied. In both
papers, we assumed that the diffusion coefficient $\sigma$ is constant,
in order to control the mixing rate of the process (see the remarks at
the end of the present paper for extensions to the nonconstant
diffusion case).

Here, we consider the estimation of the function $\sigma$, in
particular we do no more assume necessarily that it is a constant. We
observe the process in a high resolution grid, that is, $Z_{ph_n},
p=1,\ldots,n$ with $h_n \displaystyle{\mathop{\longrightarrow}_{n \rightarrow+ \infty}} 0$. As for
our previous works, we consider the case where only the position
coordinates $X_{ph_n}$ are observed (partial observations). This is of
course the main technical difficulty. This situation leads us to define
the estimator using the second- order increments of process $\Delta
_2X(p,n)=X_{(p+1)h_n}-2X_{ph_n}+X_{(p-1)h_n}$. This fact introduces
some technicalities in the proof of each result.

In the first part of the article, we consider the case of infill
statistics $t=nh_n$ is fixed. Two situations are in view: first,
$\sigma
$ is a constant and we estimate $\sigma^2$ by using a normalization of
$\Delta_2X(p,n)$, second,\vspace*{1pt} $\sigma$ is no more constant and we estimate
$ \int_0^t \sigma^2(X_s,Y_s)\,ds$. In both cases, we obtain a stable
limit theorem with rate $\sqrt n$ for the estimators (for the
definition of stable convergence in law see the next section).

This asymptotic convergence can be applied, for instance, for testing
the null hypothesis ${\mathcal H_0}$: the matrix $\sigma$ contains only
nonvanishing diagonal terms that is, $\sigma_{ij}=0$ for $i\neq j$.

In the second part, we study the infinite horizon estimation
$nh_n=t \displaystyle{\mathop{\longrightarrow}_{n \rightarrow+ \infty}} + \infty$. We assume that the
rate of mixing of the process $(Z_t, t \geq 0)$ is
sufficiently high. Whenever $\sigma$ is a constant, we obtain a central
limit theorem (CLT) for the estimator of $\sigma^2$ with rate $\sqrt
n$. However, in the case where $\sigma$ is not a constant we get a new
CLT but the rate now is $\sqrt{ nh_n}$ and the asymptotic variance is
the same as the one obtained for occupation time functionals.

The result in the infinite horizon can serve to test $\mathcal H_0:
\sigma(x,y)= \sigma$ against the sequence of alternatives $\mathcal
H^n_1: \sigma_n=\sigma+c_nd(x,y)$, for some sequence $c_n$ tending to
zero as $n$ tends to infinity, because of the difference in the
convergence rate under the null and under the sequence of alternatives.

Estimation with partial observations has been considered previously in
the literature. In \cite{GLO}, the case of one-dimensional diffusion
$V_t$ is studied. One observes only $S_t=\int_0^tV_s\,ds$, in a discrete
uniform grid. The estimation is made for the parameters defining the
variance and the drift. More recently for the same type of models, the
problem of estimation was considered in \cite{CGCR}. In this last work,
the study is nonparametric in nature; it deals with adaptive
estimation, evaluating the quadratic risk. The models in both these
articles, contrary to models of type \eqref{WuGirsanov}, do not allow
the second equation to depend on the first coordinate. It can be
written as
\[
\cases{
dS_t =  V_t \,dt, \cr\vspace*{2pt}
dV_t  =  \sigma(V_t)   \,dW_t+b(V_t)\,dt.}
\]

The literature concerning the estimation for models of type (\ref{WuGirsanov})
is rather scarce. However, two papers must be cited.
First, \cite{PSW} consider parameter estimation by using approximate
likelihoods. The horizon of estimation is infinite and they assume
$h_n \displaystyle{\mathop{\longrightarrow}_{n \rightarrow+ \infty}} 0$ and
$nh_n\displaystyle{\mathop{\longrightarrow}_{n
\rightarrow+ \infty}} + \infty$. Second, \cite{samthi} introduce, in
the case of partial observations, an Euler contrast defined using the
second coordinate only. We should point out that the present work,
while dealing with nonparametric estimation, has a nonempty
intersection with the one of \cite{samthi} when the diffusion term is
constant.

In Section~\ref{Langevin}, we consider Langevin dynamics described by
%
%e1.2 #&#
\begin{equation}
\label{LangevinEquations}
\cases{
dX_t = Y_t
\,dt,
\vspace*{4pt}\cr
dY_t  = \sqrt{2 \beta^{-1}} s(X_t)
\,dW_t-\bigl(s(X_t)s^*(X_t)Y_t+
\nabla V(X_t)\bigr) \,dt. }
\end{equation}
This form of hypo-elliptic diffusion is a particular case of (\ref{WuGirsanov}) with
$\sigma(x,y)\times\sigma^*(x,y)=\frac{2}{\beta}c(x,y)$.
This last relation is called fluctuation-dissipation relation since it
relates the magnitude of the dissipative term $- c(X_t,Y_t)   \,dt =
-s(X_t)s^*(X_t)Y_t   \,dt$ and the magnitude of
the random term $\sigma(X_t,Y_t)   \,dW_t= \sqrt{2 \beta^{-1}} s(X_t)
  \,dW_t$. The precise balance between the drift term which removes energy
in average and the stochastic term provided by the
fluctuation-dissipation relation insures that the
canonical measure is preserved by the dynamics. More precisely, under
assumptions $\mathcal{H}_0$ and $\mathcal{H}_1$ of Section~\ref{sectools}, it is proved that the solution of (\ref{LangevinEquations})
is ergodic \vspace*{2pt}with invariant probability measure
proportional to the Boltzmann \vspace*{1pt}distribution $\exp(-\beta H(x,y))$, where
$H(x,y)=\frac{1}2|y|^2+V(x)$ and $\beta$ is inversely proportional
to the temperature (see, e.g., \cite{Matti1,Tony}).

Numerical experiments are provided in Section~\ref{numexp}.

Let us end this \hyperref[sec1]{Introduction} with some comments about some possible
generalizations. In the first place, the methods that we use in this
work can be adapted for considering the power variation type estimators
defined as
\[
V_F(n)=\sum_{p=0}^{(([{t}/({2h_n})]-1)h_n}F\bigl(
\Delta_2X(p,n)\bigr),
\]
for $F$ a smooth function, usually $F(x)=|x|^r$ the $r$th
power variation (see, e.g., \cite{Jacod}). Second, it would be
possible to study an estimator constructed through a Fourier transform
method as the one defined in \cite{MM}.

%s2 #&#
\section{Tools}\label{sectools}
%s2.1 #&#
\subsection{Stable convergence}
In this article, the type of convergence we consider is the stable
convergence, introduced by Renyi, whose definition is recalled below
(see Definition~\ref{SCVdefi}).
In this subsection, all random variables or processes are defined on
some probability space
$ (\Omega, \mathbb{F},\mathbb{P} )$.

%de2.1 #&#
\begin{Definition}[(Definition~2.1 in \cite{PV})]\label{SCVdefi}
Let $Y_n$ be a sequence of random variables with values in a Polish space
$(E, \mathcal{E})$. We say that $Y_n$ converges stably with limit $Y$,
written $Y_n \displaystyle{\mathop{\longrightarrow}_{n \rightarrow+ \infty}^{\mathcal{S}}} Y$,
where $Y$ is defined on an extension $ (\Omega', \mathbb{F}',
\mathbb{P}')$ iff for any bounded, continuous function $g$ and
any bounded
$\mathbb{F}$-measurable random variable $Z$ it holds that
\[
\mathbb{E}\bigl(g(Y_n)Z\bigr) \rightarrow\mathbb{E'}
\bigl(g(Y )Z\bigr)
\]
as $n \rightarrow+ \infty$, where $\mathbb{E}$ (resp.,
$\mathbb
{E'}$) denotes the expectation with respect to probability $\mathbb{P}$
(resp., $\mathbb{P'}$).
\end{Definition}

If $\mathbb F$ is the $\sigma$-algebra generated by some random
variable (or process) $X$, then it is enough to consider $Z=h(X)$ for
some continuous and bounded $h$. It is thus clear that the stable
convergence in this situation is equivalent to the convergence in
distribution of the sequence $(Y_n,X)$ to $(Y,X)$. It is also clear
that convergence in probability implies stable convergence. As shown in
\cite{PV}, the converse holds true if $Y$ is $\mathbb F$
measurable.

Notice that we may replace the assumption $Z$ is bounded
by $Z \in\mathbb L^1(\mathbb P)$. This remark allows us to replace
$\mathbb P$ by any $\mathbb Q$ which is absolutely continuous with
respect to $\mathbb P$, that is, the following proposition.

%pr2.2 #&#
\begin{Proposition}\label{heredite}
Assume\vspace*{1pt} that $Y_n$ [defined on $ (\Omega, \mathbb{F},\mathbb
{P}
)$] converges stably to~$Y$. Let $\mathbb Q$ be a probability measure
on $\Omega$ such that $\frac{d\mathbb Q}{d \mathbb P} = H$. Then $Y_n$
[defined on $(\Omega, \mathbb{F},\mathbb{Q} )$] converges
stably to the same $Y$ [defined on $ (\Omega', \mathbb
{F'},\mathbb
{Q'}=H \mathbb P'  )$].
\end{Proposition}

In particular, in the framework of our diffusion processes, this
proposition combined with Girsanov transform theory will allow us to
``kill'' the drift.

%s2.2 #&#
\subsection{About the s.d.e. \texorpdfstring{\protect\eqref{WuGirsanov}}{(1.1)}}
In all of the paper, we will assume (at least) that the coefficients in
\eqref{WuGirsanov} satisfy:
\begin{itemize}
\item $\mathcal{H}_0$ The diffusion matrix $\sigma$ is symmetric,
smooth, bounded as well as its first and second partial derivatives and
uniformly elliptic, that is, $\forall x,y$, $\sigma(x,y) \geq\sigma_0
  \mathrm{Id}$ (in the sense of quadratic forms) for a positive constant
$\sigma
_0 >0$.
\item $\mathcal{H}_1$ The potential $V$ is lower bounded and
continuously differentiable on $\mathbb{R}^{d}$.
\item $\mathcal{H}_2$ The damping matrix $c$ is continuously
differentiable and for all $N>0$: $\sup_{|x|\leq N,   y \in
\mathbb
{R}^d} |c(x,y)| < + \infty\mbox{ and } \exists\vspace*{1pt}  c_0,   L>0  \
  c^s(x,y) \geq c_0 \mathrm{Id} $ for all $|x| >L,   y \in\mathbb{R}^d$,
$c^s$ being the symmetrization of the matrix $c$.
\end{itemize}
Under these assumptions, equation \eqref{WuGirsanov} admits an unique
strong solution which is nonexplosive. In addition, we have the
following lemma.

%le2.3 #&#
\begin{Lemma}[(Lemma~1.1 in \protect\cite{Wu})]\label{WuGigi}
Assume $\mathcal{H}_0$, $\mathcal{H}_1$ and $\mathcal{H}_2$. Then, for
every initial state $z = (x, y) \in\mathbb{R}^{2d}$, the s.d.e.
\protect\eqref{WuGirsanov} admits a unique
strong solution $\mathbb{P}_z$ (a~probability measure on $\Omega$),
which is nonexplosive. Moreover,
$\mathbb{P}_z  \ll  \mathbb{P}_z^0$ on $(\Omega, \mathcal{F}_t)$ for each
$t>0$, where $\mathbb{P}_z^0$ is the law of the solution of \eqref
{WuGirsanov} associated to $c(x, y) = 0$ and $V = 0$, and with $
(\mathcal{F}_t:=\sigma(Z_s,   0 \leq s \leq t) )_{t \geq0}$.
\end{Lemma}
%

%re2.4 #&#
\begin{remark}\label{remsig}
The formulation of $\mathcal H_0$ can be surprising. Let $\sigma^*$
denote the transposed matrix of $\sigma$. Actually the law of the
process depends on $\sigma\sigma^*$ (which is the second-order term of
the infinitesimal generator). If this symmetric matrix is smooth, it is
well known that one can find a smooth symmetric square root of it,
which is the choice we make for $\sigma$. As it will be clear in the
sequel, our estimators are related to $\sigma\sigma^*$ (hence here
$\sigma^2$).
\end{remark}

%s3 #&#
\section{Finite horizon (infill) estimation}

We consider infill estimation, that is we observe the process on a
finite time interval $[0,T]$, with a discretization step equal to $h_n$
with $h_n \displaystyle{\mathop{\longrightarrow}_{n \rightarrow+ \infty}} 0$.

According to Lemma~\ref{WuGigi} and Proposition~\ref{heredite},
\textit{any $\mathbb{P}_z^0$ stably converging sequence $Y_n$ is also
$\mathbb{P}_z$ stably converging to the same limit}. Hence, \textit{in
all of this section we will assume that $\mathcal{H}_0$ is satisfied
and that $c$ and $V$ are identically $0$.} Any result obtained in this
situation is thus true as soon as $\mathcal{H}_0$, $\mathcal{H}_1$ and
$\mathcal{H}_2$ are satisfied.

%s3.1 #&#
\subsection{The case of a constant diffusion matrix}\label{subsecons}
  We start with the definition of the ``double'' increment of the process.

Define for $1 \leq p \leq[\frac{T}{2h_n}]-1:=p_n$ (here $[\cdot]$ denotes
the integer part)
%
%e3.1 #&#
\begin{equation}
\label{acdouble}
\Delta_2X(p,n)=X_{(2p+1)h_n}-2X_{2ph_n}+X_{(2p-1)h_n}.
\end{equation}
Then
\begin{eqnarray*}
\sigma^{-1} \Delta_2X(p,n)
&= & \int_{2ph_n}^{(2p+1)h_n}W_s\,ds-\int
_{(2p-1)h_n}^{2ph_n}W_u\,du
\nonumber
\\
& = &
\int_{2ph_n}^{(2p+1)h_n}(W_s-W_{2ph_n})
\,ds+\int_{(2p-1)h_n}^{2ph_n}(W_{2ph_n}-W_u)
\,du.
\end{eqnarray*}
The right-hand side is the sum of two independent centered normal
random vectors, whose coordinates are independent, so that $\sqrt
{\frac
{3}{2 h_n^3}}   \sigma^{-1}   \Delta_2X(p,n)$ is a centered Gaussian
random vector with covariance matrix equal to Identity (recall that we
assume that $\sigma=\sigma^*$).

Furthermore, all the $(\Delta_2X(p,n))_{1\leq p \leq p_n}$
are independent (thanks to our choice of the increments).

So we define our estimator $\hat{\sigma}^2_n$ of the matrix $\sigma
^2$ as
%
%e3.2 #&#
\begin{equation}
\label{estim}
\hat{\sigma}_n^2= \frac{1}{[{T}/({2h_n})]-1}
\frac{3}{2
h_n^3} \sum_{p=1}^{[{T}/({2h_n})]-1}
\Delta_2X(p,n)\otimes\Delta _2X(p,n),
\end{equation}
where $A\otimes B$ denotes the $(d,d)$ matrix obtained by taking the
matrix product of the $(d,1)$ vector $A$ with the transposed of the
$(d,1)$ vector $B$, denoted by $B^*$.

Using what precedes, we see
that
\[
\sigma^{-1} \hat{\sigma}_n^2
\sigma^{-1}= \frac{1}{[{T}/({2h_n})]-1} \sum_{p=1}^{[{T}/({2h_n})]-1}
M(p,n),
\]
where for each $n$ the $M(p,n)$ are i.i.d. symmetric random matrices
whose entries $M_{i,j}$ are all independent for $i \geq j$, satisfying
$\mathbb{E}_z^0(M_{i,j})=\delta_{i,j}$ and $\operatorname
{Var}_z^0(M_{i,j})= 1+\delta_{i,j}$.

According to the law of large numbers and the central limit theorem for
triangular arrays of independent variables, we have the following.

%le3.1 #&#
\begin{Lemma}[(Convergence)]\label{cvas}
Assume $c=0$, $V=0$ and $\mathcal{H}_0$.
Then if\break  $h_n \displaystyle{\mathop{\longrightarrow}_{n \rightarrow+ \infty}}0$,
starting from
any initial point $z=(x,y) \in\mathbb{R}^{2d}$, we have
\[
\hat{\sigma}_n^2 \mathop{\longrightarrow}_{n \rightarrow+ \infty}^{\mathbb P_z^0}
\sigma^2,
\]
and
\[
\biggl( \biggl[\frac{T}{2h_n} \biggr]-1 \biggr)^{1/2} \bigl(
\sigma^{-1} \hat {\sigma}_n^2
\sigma^{-1} - \mathrm{Id}\bigr) \mathop{\longrightarrow}_{n \rightarrow+\infty}^{\mathcal D}
\mathcal N_{(d,d)},
\]
where $\mathcal N_{(d,d)}$ is a $(d,d)$ symmetric random matrix whose
entries are centered Gaussian random variables with $\operatorname
{Var}(\mathcal
N_{i,j})=1+\delta_{i,j}$, all the entries $(i,j)$ for $i\geq j$ being
independent.
\end{Lemma}

The consistence result is interesting since convergence in $\mathbb
P_z^0$ probability implies convergence in $\mathbb P_z$ probability
(i.e., for general $c$ and $V$). The convergence in distribution
however is not sufficient and has to be reinforced into a stable convergence.

This is the aim of what follows.

To this end, we define the sequence of processes defined for $0 \leq t
\leq T$,
%
%e3.3 #&#
\begin{equation}
\label{estimt}
\hat{\sigma}_n^2(t)= \frac{1}{[{T}/({2h_n})]-1}
\frac{3}{2 h_n^3} \sum_{p=1}^{[{t}/({2h_n})]-1}
\Delta_2X(p,n)\otimes \Delta _2X(p,n),
\end{equation}
where the empty sums are set equal to zero. We will prove
the following.

%th3.2 #&#
\begin{Theorem}[(Convergence in Skorohod's metric)]\label{cvskorohod}
Assume $c=0$, $V=0$ and $\mathcal{H}_0$.

Then if $h_n \displaystyle{\mathop{\longrightarrow}_{n \rightarrow+ \infty}}0$, starting from
any initial point $z=(x,y) \in\mathbb{R}^{2d}$, we have
\[
\biggl(W_{\cdot},\sqrt{\frac{T}{2h_n}} \bigl( \sigma^{-1}
\hat {\sigma }_{n}^2(\cdot) \sigma^{-1} - \mathrm{Id}
\bigr) \biggr) \mathop{\hbox to 2cm{\rightarrowfill}}_{n\rightarrow+ \infty}^{\mathcal{D} ([0,T] ) \times\mathcal
{D} ([0,T] ) }
(W_{\cdot},\tilde{W}_{\cdot
} ),
\]
where $(\tilde{W}_{t},   t \in[0,T] )$ is a
$(d,d)$ symmetric matrix valued random process whose entries are Wiener
processes with variance $\operatorname{Var}_{i,j}(t)=(1+\delta_{i,j})
  (t/T)$, all
the entries $(i,j)$ for $i\geq j$ being independent. In addition
${\tilde{W}}_{.}$ is independent of $W_{.}$.
\end{Theorem}

According to the discussion on stable convergence, we immediately
deduce the following.

%co3.3 #&#
\begin{Corollary}[(Stable convergence)]\label{theoconstantcVnuls}
Under assumptions $\mathcal{H}_0$, $\mathcal{H}_1$ and
$\mathcal{H}_2$, if $h_n \displaystyle{\mathop{\longrightarrow}_{n \rightarrow+
\infty}}0$, starting from any initial point $z=(x,y) \in\mathbb{R}^{2d}$,
we have
\[
\sqrt{\frac{T}{2h_n}} \bigl(\sigma^{-1} \hat{\sigma}_n^2
\sigma ^{-1} - \mathrm{Id} \bigr) \mathop{\longrightarrow}_{n \rightarrow+ \infty}^{\mathcal{S}}
\mathcal{N}_{(d,d)},
\]
where $\mathcal{N}_{(d,d)}$ is as in Lemma~\ref{cvas}.
\end{Corollary}

\begin{pf*}{Proof of Theorem~\protect\ref{cvskorohod}}
In the following, we fix $T=1$ without loss of generality. Notice that
we may also replace $\frac{T}{2h_n}$ by $[\frac{T}{2h_n}]-1$ (using
Slutsky's theorem if one wants).

The convergence of
\[
t \mapsto Z_n(t)=\sqrt{\frac{1}{2h_n}} \sigma^{-1}
\hat {\sigma }_n^2(t) \sigma^{-1}
\]
to a matrix of Wiener processes is proved as for Donsker invariance
principle. The only difference here is that instead of an i.i.d. sample
we look at a triangular array of i.i.d. random vectors (on each row),
but the proof in \cite{Bil} applies in this situation. This result is
sometimes called Donsker--Prohorov invariance principle. Writing $W_t$
as the sum of its increments on the grid given by the intervals
$[(2p-1)h_n,(2p+1)h_n]$ the convergence of the joint law of
$(W_.,Z_n(\cdot))$ in $\mathcal D([0,1])$ is proved in exactly the same
way.

The final independence assumption is a simple covariance
calculation.
\end{pf*}

%s3.2 #&#
\subsection{Estimation of the noise, general case}

In this section, we do not assume anymore that the diffusion
term $\sigma$ is constant.

In the following, we want to estimate
$\int_0^t \sigma^2(X_s,Y_s)\,ds$, for any $0 \leq t \leq T$.

To this
end, we introduce the quadratic variation process defined for $n \in
\mathbb{N}^*$ and $0 \leq t \leq T$ as
%
%e3.4 #&#
\begin{equation}
\label{variation}
\mathcal{QV}_{h_n}(t)=\frac{1}{h_n^2} \sum
_{p=1}^{ [{t}/({2h_n})]-1} \Delta_2X(p,n)\otimes
\Delta_2X(p,n),
\end{equation}
with $\Delta_2X(p,n)$ defined in (\ref{acdouble}). The main result of
this section is the following.

%th3.4 #&#
\begin{Theorem}\label{thmmainvarinfill}
Under assumptions $\mathcal{H}_0$, $\mathcal{H}_1$ and
$\mathcal{H}_2$, if $h_n \displaystyle{\mathop{\longrightarrow}_{n \rightarrow+
\infty}}0$, starting from any initial point $z=(x,y) \in\mathbb{R}^{2d}$,
we have for any $0 \leq t \leq T$
\[
\mathcal{QV}_{h_n}(t) \mathop{\longrightarrow}_{n \rightarrow+ \infty}^{\mathbb
{P}_z}
\frac{1}3 \int_0^t
\sigma^2(X_s,Y_s)\,ds,
\]
and
\[
\sqrt{\frac{1}{h_n}} \biggl(\mathcal{QV}_{h_n}(t)-
\frac{1}3 \int_0^t \sigma^2
(X_s,Y_s)\,ds \biggr)\mathop{\longrightarrow}_{n \rightarrow+
\infty}^{\mathcal S}
\frac{2}3 \int_0^t
\sigma(X_s,Y_s) \,d\tilde W_s
\sigma(X_s,Y_s),
\]
where $ ( \tilde{W}_t, t \in[0,T] )$ is a
symmetric matrix valued random process independent of the initial
Wiener process $W_{.}$, whose entries $\tilde W_.(i,j)$ are Wiener
processes with variance $V_{i,j}(t)=(1+\delta_{i,j})t$, these entries
being all independent for $i\geq j$.
\end{Theorem}

Recall that for the proof of this theorem, we only need to consider the
case where $c=0$ and $V=0$.

In this case, the strong solution, with
initial conditions $(X_0,Y_0)=(x,y)=z$, can be written as
\[
Z_t=(X_t,Y_t)= \biggl( x +yt+ \int
_0^t Y_s\,ds,y+ \int
_0^t \sigma (X_s,Y_s)
\,dW_s \biggr).
\]
We thus have
\begin{eqnarray*}
\Delta_2X(p,n)&=&\int_{2ph_n}^{(2p+1)h_n} \biggl[
\int_0^s\sigma (X_u,Y_u)
\,dW_u \biggr]\,ds\\
&&{}-\int_{(2p-1)h_n}^{2ph_n}
\biggl[\int_0^s\sigma (X_u,Y_u)
\,dW_u \biggr]\,ds.
\end{eqnarray*}

Using Fubini's theorem for stochastic integrals, one gets
\begin{eqnarray*}
\Delta_2X(p,n) & = & h_n
\int_0^{2ph_n}\sigma(X_u,Y_u)
\,dW_u\\
&&{}+\int_{2ph_n}^{(2p+1)h_n}
\bigl((2p+1)h_n-u\bigr)\sigma(X_u,Y_u)
\,dW_u
\\
&&{} -  h_n\int_0^{(2p-1)h_n}
\sigma(X_u,Y_u)\,dW_u\\
&&{}-\int
_{(2p-1)h_n}^{2ph_n}(2ph_n-u)
\sigma(X_u,Y_u)\,dW_u,
\end{eqnarray*}
thus
%
%e3.5 #&#
\begin{equation}
\label{ecriture} \Delta_2X(p,n)= \int_{(2p-1)h_n}^{(2p+1)h_n}\bigl(h_n-|u-2ph_n|\bigr)
\sigma (X_u,Y_u)\,dW_u.
\end{equation}
If $p\neq q$ are two integers, denoting by $\Delta_2X(p,n,i)$ the
$i$th coordinate of $\Delta_2X(p,n)$, we immediately have, for all
$i,j=1,\ldots,d$,
%
%e3.6 #&#
\begin{equation}
\label{covnul}
\mathbb{E}_z^0 \bigl(\Delta_2X(p,n,i)
\Delta_2X(q,n,j) \bigr) = 0.
\end{equation}

As a warm up lap, we look at the convergence of the first moment of
$\mathcal{QV}_{h_n}$.

%le3.5 #&#
\begin{Lemma}[(Preliminary result)]\label{moment1}
Assume $c=0$, $V=0$ and $\mathcal H_0$. Then, if $h_n \displaystyle{\mathop{\longrightarrow}_{n
\rightarrow+ \infty}}0$, starting from any initial point $z=(x,y)
\in
\mathbb{R}^{2d}$, we have for any $0 \leq t \leq T$,
\[
\mathbb{E}_z^0 \mathcal{QV}_{h_n}(t) \mathop{\longrightarrow}_{n\rightarrow+ \infty}\frac{1}3 \int_0^t
\mathbb{E}_z^0\sigma^2(X_u,Y_u)
\,du.
\]
\end{Lemma}

Recall that we assumed $\sigma=\sigma^*$, and of course look at the
previous equality as an equality between real matrices.

\begin{pf*}{Proof of Lemma~\protect\ref{moment1}}
First, using It\^{o}'s
isometry and equality \eqref{ecriture}, one gets
\[
\mathbb{E}_z^0 \bigl(\Delta_2X(p,h_n)
\otimes\Delta_2X(p,h_n) \bigr)=\int_{(2p-1)h_n}^{(2p+1)h_n}\bigl(h_n-|u-2ph_n|\bigr)^2
\mathbb{E}_z^0\sigma^2(X_u,Y_u)
\,du.
\]
Since
\[
\int_{(2p-1)h_n}^{(2p+1)h_n}\bigl(h_n-|u-2ph_n|\bigr)^2
\,du = \frac{2}3 h_n^3,
\]
we thus have
\begin{eqnarray*}
&& \frac{1}{h_n^2} \mathbb{E}_z^0 \bigl(
\Delta_2X(p,h_n)\otimes \Delta _2X(p,h_n)
\bigr) - \frac{1}3 \int_{(2p-1)h_n}^{(2p+1)h_n}
\mathbb{E} _z^0\sigma ^2(X_u,Y_u)
\,du\\
&&\qquad =\frac{1}{h^2_n} \int_{(2p-1)h_n}^{(2p+1)h_n}\bigl(h_n-|u-2ph_n|\bigr)^2
\\
&&\qquad\quad{}\times\mathbb{E} _z^0 \bigl(\sigma^2(X_u,Y_u)-
\sigma ^2(X_{(2p-1)h_n},Y_{(2p-1)h_n}) \bigr) \,du\\
&&\qquad\quad {}+
\frac{1}3 \int_{(2p-1)h_n}^{(2p+1)h_n}\mathbb{E}
_z^0 \bigl(\sigma^2(X_{(2p-1)h_n},Y_{(2p-1)h_n})-
\sigma ^2(X_u,Y_u) \bigr) \,du.
\end{eqnarray*}

Define on $\Omega\times[0,t]$, the sequence of random (matrices)
\[
G_{n}(u)=\sum_{p=1}^{[{t}/({2h_n})]-1}
\sigma ^2(X_{(2p-1)h_n},Y_{(2p-1)h_n}) \BBone_{(2p-1)h_n\leq u<(2p+1)h_n}.
\]

Since $\sigma$ is continuous, $G_n$ converges $\mathbb P_z^0\otimes \,du$
almost everywhere to $\sigma^2(X_u,Y_u)$. In addition, since $\sigma$
is bounded, $G_n$ is dominated by a constant which is $\mathbb
P_z^0\otimes \,du$ integrable on $\Omega\times[0,t]$.

Hence, using
Lebesgue bounded convergence theorem, we get that
\[
\int_0^t \mathbb{E}_z^0
\bigl(G_n(u) - \sigma^2(X_u,Y_u)
\bigr) \,du \to  0.
\]

Similarly, the variables
\begin{eqnarray*}
&&\sum_{p=1}^{[{t}/({2h_n})]-1} \frac{(h_n-|u-2ph_n|)^2}{h_n^2}
\\
&&\qquad{}\times\BBone_{(2p-1)h_n\leq u<(2p+1)h_n} \bigl(\sigma ^2(X_{(2p-1)h_n},Y_{(2p-1)h_n})-
\sigma^2(X_u,Y_u)\bigr)
\end{eqnarray*}
are bounded and converge almost everywhere to $0$, so that their
expectation also goes to $0$. This completes the proof.
\end{pf*}

Of course, a careful look at this proof shows that we did not use all
the strength of $\mathcal H_0$, only the fact that $\sigma$ is
continuous and bounded. It is thus clearly possible to improve upon
this result, using the same idea of introducing the skeleton Markov
chain and controlling the errors.

Hence, introduce
%
%e3.7 #&#
\begin{equation}\qquad
\label{eqH} \Delta_2 H(p,h_n)= \int
_{(2p-1)h_n}^{(2p+1)h_n} \bigl(h_n-|u-2ph_n|\bigr)
\sigma(X_{(2p-1)h_n},Y_{(2p-1)h_n}) \,dW_u.\hspace*{-6pt}
\end{equation}
We may decompose
%
%e3.8 #&#
\begin{eqnarray}\label{eqdecomp}
J^1_n + J^2_n +
J^3_n &=& \mathcal{QV}_{h_n}(t)-
\frac{1}3 \int_0^t
\sigma^2(X_u,Y_u)\,du
\end{eqnarray}
with
\begin{eqnarray*}
J^1_n &=& \mathcal{QV}_{h_n}(t) -
\frac{1}{h_n^2} \sum_{p=1}^{
[{t}/({2h_n})]-1}
\Delta_2 H(p,h_n)\otimes\Delta_2
H(p,h_n),
\nonumber
\\
\nonumber
J^2_n &=& \Biggl(\frac{1}{h_n^2} \sum
_{p=1}^{ [{t}/({2h_n}) ]-1} \Delta_2 H(p,h_n)
\otimes\Delta_2 H(p,h_n) \Biggr) - \frac{1}3
\biggl(\int_0^t G_n(u)\,du \biggr),
\\
J^3_n &=& \frac{1}3 \biggl(\int
_0^t G_n(u)\,du - \int
_0^t \sigma ^2(X_u,Y_u)
\,du \biggr).
\nonumber
\end{eqnarray*}
For $A= (A_{i,j} )_{1 \leq i \leq q, 1 \leq j \leq r}$ a
$q\times r$ real matrix, we define $|A|$ as $|A|= \break \max_{1 \leq i \leq q,
1 \leq j \leq r}|A_{i,j}|$.\vadjust{\goodbreak}

We then have the following.

%le3.6 #&#
\begin{Lemma}\label{lemapprox}
Assume $c=0$, $V=0$, $(X_0,Y_0)=(x,y) \in\mathbb{R}^{2d}$ and
$\mathcal H_0$. Then there exist constants $C$ depending on $\sigma$,
its derivatives and the dimension only, such that for any $0 \leq t
\leq T$,
%
%e3.9 #&#
\begin{equation}
\label{eqapprox1}
\mathbb{E}_z^0 \biggl(\biggl\llvert \int
_0^t G_n(u)\,du - \int
_0^t \sigma ^2(X_u,Y_u)
\,du\biggr\rrvert \biggr) \leq Ct \sqrt{h_n},
\end{equation}
and
%
%e3.10 #&#
\begin{equation}
\label{eqapprox2}
\mathbb{E}_z^0 \bigl(\bigl\llvert
\Delta_2 X(p,h_n)-\Delta_2
H(p,h_n)\bigr\rrvert ^2 \bigr) \leq C
h_n^4.
\end{equation}
\end{Lemma}
\begin{pf}
For the first part, it is enough to show that
\[
\int_0^t \mathbb{E}_z^0
\bigl\llvert G_n(u) - \sigma^2(X_u,Y_u)
\bigr\rrvert \,du \leq Ct \sqrt{h_n}.
\]
But using the fact that $\sigma$ and its first derivatives are
continuous and bounded, there exists a constant $C$ only depending on
these quantities (but which may change from line to line), such that
%
%e3.11 #&#
\begin{equation}
\label{eqcontinu}
\qquad\int_0^t \mathbb{E}_z^0
\bigl\llvert G_n(u) - \sigma^2(X_u,Y_u)
\bigr\rrvert \,du \leq Ct \sup_{|a-b|\leq2h_n} \mathbb{E}_z^0\bigl(|Z_a-Z_b|\bigr)
\leq Ct \sqrt{h_n}.\hspace*{-12pt}
\end{equation}
For the second part, we have
\begin{eqnarray*}
&& \mathbb{E}_z^0 \bigl(\bigl\llvert \Delta_2 X(p,h_n)-\Delta_2
H(p,h_n)\bigr\rrvert ^2 \bigr)
\\
&&\qquad= \mathbb{E}_z^0 \biggl(\int_{(2p-1)h_n}^{(2p+1)h_n}
\bigl(h_n-|u-2ph_n|\bigr)^2 \\
&&\qquad\quad{}\times\operatorname{Trace}\bigl(\bigl(
\sigma(X_{(2p-1)h_n},Y_{(2p-1)h_n})-\sigma(X_u,Y_u)
\bigr)^2\bigr) \,du \biggr),
\end{eqnarray*}
from which the result easily follows as before.
\end{pf}
We deduce immediately
the following.

%pr3.7 #&#
\begin{Proposition}\label{propapprox}
Assume $c=0$, $V=0$ and $\mathcal H_0$. Then, if $h_n
\displaystyle{\mathop{\longrightarrow}_{n \rightarrow+ \infty}}0$, starting from any initial point $z=(x,y)
\in
\mathbb{R}^{2d}$, we have for any\vspace*{1.5pt} $0 \leq t \leq T$, that $J_n^1$ and
$J_n^3$ are converging to $0$ in $\mathbb L^1(\mathbb P_z^0)$ (with
rates $h_n$ and $\sqrt{h_n}$), hence in $\mathbb P_z^0$ probability.
\end{Proposition}

\begin{pf}
The result for $J_n^3$ is contained in the previous lemma. For $J_n^1$,
we calculate $\mathbb{E}_z^0[|J_n^1|]$. The $(i,j)$th term of
$J_n^1$ is
given by
\[
\frac{1}{h_n^2} \sum_{p=1}^{ [{t}/({2h_n}) ]-1} \bigl(
\Delta _2 X(p,h_n,i)-\Delta_2
H(p,h_n,i)\bigr) \bigl(\Delta_2 X(p,h_n,j)-
\Delta_2 H(p,h_n,j)\bigr),
\]
so that, according to the previous lemma and the Cauchy--Schwarz
inequality, we thus have $\mathbb{E}_z^0[|J_n^1|] \leq Ct   h_n$.
\end{pf}

In order to prove the first part of Theorem~\ref{thmmainvarinfill},
that is, the convergence in probability, it remains to look at $J_n^2$.
We have
%
%e3.12 #&#
\begin{eqnarray}
J_n^2 &=& \frac{1}{h_n^2} \sum
_{p=1}^{ [{t}/({2h_n})
]-1} \sigma(X_{(2p-1)h_n},Y_{(2p-1)h_n})
\biggl( M(p,h_n) - \frac{2h_n^3}{3} \mathrm{Id} \biggr)
\nonumber
\\[-8pt]
\label{jn2}
\\[-8pt]
\nonumber
&&{}\times\sigma(X_{(2p-1)h_n},Y_{(2p-1)h_n}),
\end{eqnarray}
where
\[
M(p,h_n) = \Delta_2 W(p,h_n)\otimes
\Delta_2 W(p,h_n)
\]
and
\[
\Delta_2 W(p,h_n) = \int_{(2p-1)h_n}^{(2p+1)h_n}
\bigl(h_n-|u-2ph_n|\bigr) \,dW_u.
\]
As before, we start with an estimation lemma.

%le3.8 #&#
\begin{Lemma}\label{lemapprox2}
Assume $c=0$, $V=0$, $(X_0,Y_0)=(x,y) \in\mathbb{R}^{2d}$ and
$\mathcal H_0$. Then there exist constants $C$ depending on $\sigma$,
its derivatives and the dimension only, such that
%
%e3.13 #&#
\begin{equation}
\label{eqapprox3}
\mathbb{E}_z^0 \biggl(\biggl\llvert
M(p,h_n) - \frac{2h^3_n}{3} \mathrm{Id}\biggr\rrvert ^2 \biggr)
\leq C h^6_n.
\end{equation}
\end{Lemma}

\begin{pf}
We shall look separately at the diagonal terms and the off diagonal
terms of $M(p,h_n)   -   \frac{2h^3_n}{3}   \mathrm{Id}$.

The off diagonal
terms are of the form $A_{i,j}(n)$ with
%
%e3.14 #&#
\begin{eqnarray}
A_{i,j}(n)&=& \biggl(\int_{(2p-1)h_n}^{(2p+1)h_n}
\bigl(h_n-|u-2ph_n|\bigr) \,dW^i_u
\biggr)
\nonumber
\\[-8pt]
\label{offdiag}
\\[-8pt]
\nonumber
&&{}\times \biggl(\int_{(2p-1)h_n}^{(2p+1)h_n} \bigl(h_n-|u-2ph_n|\bigr)
\,dW^j_u \biggr),
\end{eqnarray}
where $W^i$ and $W^j$ are independent linear Brownian motions.
Introduce the martingales
\[
U_i(s)=\int_{(2p-1)h_n}^{s}
\bigl(h_n-|u-2ph_n|\bigr) \,dW^i_u
\]
defined for $(2p-1)h_n \leq s \leq(2p+1)h_n$. Using It\^{o}'s formula,
$A_{i,j}(n)$ can be rewritten
\begin{eqnarray*}
&& \biggl(\int_{(2p-1)h_n}^{(2p+1)h_n} \bigl(U_j(u)-U_j
\bigl((2p-1)h_n\bigr)\bigr) \bigl(h_n-|u-2ph_n|\bigr)
\,dW^i_u \biggr)
\\
&&\qquad{}+ \biggl(\int_{(2p-1)h_n}^{(2p+1)h_n} \bigl(U_i(u)-U_i
\bigl((2p-1)h_n\bigr)\bigr) \bigl(h_n-|u-2ph_n|\bigr)
\,dW^j_u \biggr)
\end{eqnarray*}
so that
\begin{eqnarray*}
&&\mathbb{E}_z^0 \bigl(A^2_{i,j}(n)
\bigr)\\
 &&\qquad= 2 \int_{(2p-1)h_n}^{(2p+1)h_n} \bigl(h_n-|u-2ph_n|\bigr)^2
\mathbb{E}_z^0\bigl[\bigl(U_i(u)-U_i
\bigl((2p-1)h_n\bigr)\bigr)^2\bigr] \,du
\\
&&\qquad= 2 \int_{(2p-1)h_n}^{(2p+1)h_n} \bigl(h_n-|u-2ph_n|\bigr)^2
\\
&&\qquad\quad{}\times\biggl(\int_{(2p-1)h_n}^{u} \bigl(h_n-|s-2ph_n|\bigr)^2
\,ds \biggr) \,du
\\
&&\qquad= C h_n^6,
\end{eqnarray*}
where $C$ is some universal constant, so that we get the result.

The diagonal terms can be written $A_{i,i}(n)$ with
%
%e3.15 #&#
\begin{eqnarray}
A_{i,i}(n)&=& \biggl(\int_{(2p-1)h_n}^{(2p+1)h_n}
\bigl(h_n-|u-2ph_n|\bigr) \,dW^i_u
\biggr)^2 - \frac{2}3 h_n^3
\nonumber
\\[-8pt]
\label{diag}
\\[-8pt]
&=& 2 \int_{(2p-1)h_n}^{(2p+1)h_n} \bigl(U_i(u)-U_i
\bigl((2p-1)h_n\bigr)\bigr) \bigl(h_n-|u-2ph_n|\bigr)
\,dW^i_u,
\nonumber
\end{eqnarray}
and we can conclude exactly as before.
\end{pf}

We can now state
the following.

%pr3.9 #&#
\begin{Proposition}\label{propapprox2}
Assume $c=0$, $V=0$ and $\mathcal H_0$. Then, if $h_n
\displaystyle{\mathop{\longrightarrow}_{n\rightarrow+ \infty}}0$, starting from any initial point $z=(x,y)
\in
\mathbb{R}^{2d}$, we have for any\vspace*{1pt} $0 \leq t \leq T$, that $J_n^2$
converges to $0$ in $\mathbb L^2(\mathbb P_z^0)$ (with rate $\sqrt
{h_n}$), hence in $\mathbb P_z^0$ probability.
\end{Proposition}

\begin{pf}
We look at each term $(J_n^2)_{ij}$ of the matrix $J_n^2$. Such a term
can be written in the form
\[
\bigl(J_n^2\bigr)_{ij} =\frac{1}{h_n^2}
\sum_{p=1}^{ [{t}/({2h_n})
]-1} \sum
_{l,k=1}^d a_{l,k,i,j}(X_{(2p-1)h_n},Y_{(2p-1)h_n})
A_{l,k}(p,n),
\]
where the $a_{l,k,i,j}$'s are $C^2_b$ functions, and the $A_{l,k}(p,n)$
are defined in the proof of the previous lemma (here we make explicit
the dependence in $p$). Hence,
\begin{eqnarray*}
&& h_n^4 \bigl(J_n^2
\bigr)^2_{ij} \\
&&\qquad=\sum_{p,q=1}^{ [{t}/({2h_n})]-1} \sum
_{l,k,i,j=1}^d b_{l,k,i,j}(X_{(2p-1)h_n},Y_{(2p-1)h_n})
\\
&&\qquad\quad{}\times c_{l,k,i,j}(X_{(2q-1)h_n},Y_{(2q-1)h_n}) A_{l,k}(p,n)
A_{i,j}(q,n)
\end{eqnarray*}
for some new functions $b_{l,k,i,j}$ and $c_{l,k,i,j}$. As we remarked
in \eqref{covnul}, the expectation of terms where $p\neq q$ is equal to
$0$, so that
\begin{eqnarray*}
\mathbb{E}_z^0 \bigl[h_n^4
\bigl(J_n^2\bigr)^2_{ij} \bigr] &
\leq& C \sum_{l,k,i,j=1}^d \sum
_{p=1}^{ [{t}/({2h_n})]-1} \mathbb{E} _z^0
\bigl[A_{l,k}(p,n) A_{i,j}(p,n) \bigr]
\\
&\leq& Ct h_n^5,
\end{eqnarray*}
according to the previous lemma and the Cauchy--Schwarz inequality, hence
the result.
\end{pf}

We thus have obtained the first part of the main theorem,
that is, the following.

%co3.10 #&#
\begin{Corollary}[(Consistence result)]\label{consistence}
Under assumptions $\mathcal{H}_0$, $\mathcal{H}_1$ and
$\mathcal{H}_2$, if $h_n
\displaystyle{\mathop{\longrightarrow}_{n \rightarrow+ \infty}}0$, starting from any initial point $z=(x,y) \in\mathbb{R}^{2d}$,
we have for any $0 \leq t \leq T$
\[
\mathcal{QV}_{h_n}(t) \mathop{\longrightarrow}_{n \rightarrow+ \infty}^{\mathbb
{P}_z}
\frac{1}3 \int_0^t
\sigma^2(X_s,Y_s)\,ds.
\]
\end{Corollary}

We turn now to the second part of the main theorem, that is, the
obtention of confidence intervals. Again we assume first that $c=0$ and
$V=0$.

Since\vspace*{1pt} we will normalize by $\sqrt{h_n,}$ we immediately see
that the first ``error'' term $J_n^1/\sqrt{h_n}$ converges to $0$ in
$\mathbb P_z^0$ probability according to the rate of convergence we
obtained in Proposition~\ref{propapprox}.

For the second error term
$J_n^3$, the convergence rate in $\sqrt{h_n}$ is not sufficient to
conclude. So, we have to improve on it.

%le3.11 #&#
\begin{Lemma}\label{lemJ3}
Assume $c=0$, $V=0$, $(X_0,Y_0)=(x,y) \in\mathbb{R}^{2d}$ and
$\mathcal H_0$. Then there exists some constant $C$ depending on
$\sigma$,
its first two derivatives and the dimension only, such that for any
$0 \leq t \leq T$,
%
%e3.16 #&#
\begin{equation}
\label{eqapprox11}
\mathbb{E}_z^0 \biggl(\biggl\llvert \int
_0^t G_n(u)\,du - \int
_0^t \sigma ^2(X_u,Y_u)
\,du\biggr\rrvert \biggr) \leq Ct h_n,
\end{equation}
hence $ (\int_0^t   G_n(u)\,du - \int_0^t   \sigma
^2(X_u,Y_u)\,du )/\sqrt{h_n}$ goes to $0$ in $\mathbb P_z^0$ probability.
\end{Lemma}

\begin{pf}
To begin with
\begin{eqnarray*}
&&\sigma^2(X_u,Y_u) - G_{n}(u)\\
&&\qquad=
\sum_{p=1}^{[{t}/({2h_n})]-1} \bigl(\sigma
^2(X_u,Y_u)-\sigma^2(X_{(2p-1)h_n},Y_{(2p-1)h_n})
\bigr) \BBone _{(2p-1)h_n\leq u<(2p+1)h_n}.
\end{eqnarray*}
Now look at each coordinate, and to simplify denote by $f$ the
coefficient $\sigma^2_{ij}$. It holds
\begin{eqnarray*}
&& f(Z_u) - G_{n}^{ij}(u)\\
&&\qquad=\sum_{p=1}^{[{t}/({2h_n})]-1}   \BBone_{(2p-1)h_n\leq u<(2p+1)h_n}
\int_{(2p-1)h_n}^{u} \biggl(\bigl\langle\sigma(Z_s)
\nabla_y f(Z_s), dW_s\bigr\rangle \\
&&\qquad\quad{}+ \frac{1}2
\operatorname{Trace}\bigl(\sigma D^2_y f \sigma\bigr) (Z_s)
\,ds + \bigl\langle Y_s,\nabla_x f(Z_s) \bigr\rangle\, ds \biggr)
\\
&&\qquad= \sum_{p=1}^{[{t}/({2h_n})]-1} \BBone _{(2p-1)h_n\leq
u<(2p+1)h_n}
\bigl(I^1(n,p,u)+I^2(n,p,u)\\
&&\qquad\quad{}+I^3(n,p,u)+I^4(n,p,u)
\bigr),
\end{eqnarray*}
with
\begin{eqnarray*}
I^1(n,p,u)&=& \int_{(2p-1)h_n}^{u} \bigl\langle
\sigma(Z_s) \nabla_y f(Z_s)-
\sigma(Z_{(2p-1)h_n}) \nabla_y f(Z_{(2p-1)h_n}),
dW_s\bigr\rangle,
\\
I^2(n,p,u)&=& \bigl\langle\sigma(Z_{(2p-1)h_n}) \nabla_y
f(Z_{(2p-1)h_n}), W_u-W_{(2p-1)h_n}\bigr\rangle,
\\
I^3(n,p,u)&=& \int_{(2p-1)h_n}^{u}
\frac{1}2 \operatorname{Trace}\bigl(\sigma D^2_y f \sigma\bigr)
(Z_s) \,ds,
\\
I^4(n,p,u)&=& \int_{(2p-1)h_n}^{u}
 \bigl\langle Y_s,\nabla_x f(Z_s) \bigr\rangle\, ds.
\end{eqnarray*}

Notice that $|I^3(n,p,u)|\leq C (u-(2p-1)h_n)$ so that
\[
\int_0^t \sum_{p=1}^{[{t}/({2h_n})]-1}
\BBone _{(2p-1)h_n\leq
u<(2p+1)h_n} \bigl|I^3(n,p,u)\bigr| \,du \leq Ct h_n.
\]

Similarly, $|I^4(n,p,u)| \leq C   (\sup_{0\leq s \leq t} |Y_s|)
(u-(2p-1)h_n)$ so that
\begin{eqnarray*}
&&\mathbb{E}_z^0 \Biggl(\int_0^t
\sum_{p=1}^{[{t}/({2h_n})]-1} \BBone _{(2p-1)h_n\leq u<(2p+1)h_n}
\bigl|I^4(n,p,u)\bigr| \,du \Biggr)\\
&&\qquad \leq Ct h_n
\mathbb{E}_z^0 \Bigl(\sup_{0\leq s \leq t}
|Y_s| \Bigr) \leq Ct\bigl(1+t^{1/2}\bigr) h_n
\end{eqnarray*}
according\vadjust{\goodbreak} to the Burkholder--Davis--Gundy inequality.

Now
\begin{eqnarray*}
&&\bigl(\mathbb{E}_z^0\bigl(\bigl|I^1(n,p,u)\bigr|
\bigr)\bigr)^2\\
 &&\qquad\leq E_z^0
\bigl(\bigl|I^1(n,p,u)\bigr|^2\bigr)
\\
&&\qquad= E_z^0 \biggl[\int_{(2p-1)h_n}^{u}
\bigl|\sigma(Z_s) \nabla_y f(Z_s)-\sigma
(Z_{(2p-1)h_n}) \nabla_y f(Z_{(2p-1)h_n}\bigr|^2 \,ds
\biggr]
\\
&&\qquad\leq C \bigl(u-(2p-1)h_n\bigr) \mathbb{E}_z^0
\Bigl(\sup_{|a-b|\leq2h_n} |Z_a-Z_b|^2
\Bigr)
\\
&&\qquad\leq C h_n \bigl(u-(2p-1)h_n\bigr)
\end{eqnarray*}
using the fact that $\sigma$ and its first two derivatives are bounded
and \eqref{eqcontinu}. It follows that
\[
\mathbb{E}_z^0 \Biggl(\int_0^t
\sum_{p=1}^{[{t}/({2h_n})]-1} \BBone _{(2p-1)h_n\leq u<(2p+1)h_n}
\bigl|I^1(n,p,u)\bigr| \,du \Biggr) \leq Ct h_n.
\]

Finally,
\begin{eqnarray*}
&& \Biggl(\mathbb{E}_z^0\Biggl\llvert \int_0^t   \sum_{p=1}^{[{t}/({2h_n})]-1}
\BBone
_{(2p-1)h_n\leq u<(2p+1)h_n}   I^2(n,p,u)   \,du \Biggr\rrvert  \Biggr)^2 \\
&&\qquad \leq
\mathbb{E}_z^0 \Biggl(\Biggl\llvert \int_0^t   \sum_{p=1}^{[{t}/({2h_n})]-1}
\BBone_{(2p-1)h_n\leq u<(2p+1)h_n}   I^2(n,p,u)   \,du\Biggr\rrvert ^2 \Biggr)
\\
&&\qquad\leq 2 \mathbb{E}_z^0 \Biggl[\int_0^t\!
\int_0^t \sum_{p,q=1}^{[{t}/({2h_n})]-1}
\BBone_{(2p-1)h_n\leq u<(2p+1)h_n}\\
&&\qquad\quad{}\times \BBone_{(2q-1)h_n\leq
s<(2q+1)h_n} \BBone_{s\leq u}
I^2(n,p,u) I^2(n,q,s) \,ds \,du \Biggr].
\end{eqnarray*}
As before, if $(2p-1)h_n\leq u<(2p+1)h_n$ and $(2q-1)h_n\leq
s<(2q+1)h_n$,
\[
\mathbb{E}_z^0\bigl(I^2(n,p,u)
I^2(n,q,s)\bigr)=0
\]
as soon as $p\neq q$.

If $p=q$,
\[
\bigl|\mathbb{E}_z^0\bigl(I^2(n,p,u)
I^2(n,p,s)\bigr)\bigr|\leq C \sqrt {u-(2p-1)h_n} \sqrt
{s-(2p-1)h_n},
\]
so that for a fixed $u$ between $(2p-1)h_n$ and $(2p+1)h_n$, $s$
belongs to $[(2p-1)h_n,u]$ and
\[
\int_{(2p-1)h_n}^u \bigl|\mathbb{E}_z^0
\bigl(I^2(n,p,u) I^2(n,p,s)\bigr)\bigr| \,ds \leq C
h_n^{3/2} \bigl(u-(2p-1)h_n\bigr)^{1/2}.
\]
Integrating with respect to $du$, we finally get
\[
\Biggl(\mathbb{E}_z^0\Biggl\llvert \int
_0^t \sum_{p=1}^{[{t}/({2h_n})]-1}
\BBone _{(2p-1)h_n\leq u<(2p+1)h_n} I^2(n,p,u) \,du\Biggr\rrvert
\Biggr)^2 \leq Ct h_n^2,
\]
as expected.
\end{pf}

We turn now to the central limit theorem for $J_n^2$ defined in \eqref
{jn2}. We will prove
the following.

%pr3.12 #&#
\begin{Proposition}\label{clt}
Assume $c(x,y)=0$, $V=0$ and $\mathcal{H}_0$. If $h_n \displaystyle{\mathop{\longrightarrow}_{n\rightarrow+ \infty}}0$,
starting from any initial point $z=(x,y)
\in
\mathbb{R}^{2d}$ and $\forall  0 \leq t \leq T$,
\[
\sqrt{\frac{1}{h_n}} J_n^2(t) \mathop{
\longrightarrow}_{n \rightarrow+ \infty}^{\mathcal S} \frac{2}3 \int
_0^t \sigma (X_u,Y_u)
\,d\tilde W_u \sigma(X_u,Y_u),
\]
where $(\tilde{W}_t,   t \in[0,T] )$ is a
symmetric matrix valued random process independent of the initial
Wiener process $W_{.}$, whose entries $\tilde W_.(i,j)$ are Wiener
processes with variance $V_{i,j}(t)=(1+\delta_{i,j})t$, these entries
being all independent for $i\geq j$.
\end{Proposition}

\begin{pf}
Define
\[
\xi_{n,p}= \frac{1}{h_n^2} \sigma(X_{(2p-1)h_n},Y_{(2p-1)h_n})
\biggl( M(p,h_n) - \frac{2h_n^3}{3} \mathrm{Id} \biggr) \sigma
(X_{(2p-1)h_n},Y_{(2p-1)h_n}),
\]
and $\mathcal G_{n,p}$ the $\sigma$-field generated by the $\xi_{n,j}$
for $j\leq p$. As we already saw
\[
\mathbb{E}_z^0 [\xi_{n,p}|\mathcal
G_{n,p-1} ] = 0
\]
(here the null matrix), saying that for a fixed $n$ the $\xi_{n,p}$ are
martingale increments and $J^2_n(t) =\sum_p   \xi_{n,p}$.

In order to prove the proposition, we can first show that for
all $N \in\mathbb N$, all $N$-uple $t_1 <\cdots,t_N \leq t$,
\[
\sqrt{\frac{1}{h_n}} \bigl(J_n^2(t),
W_{t_1}, \ldots, W_{t_N} \bigr) \mathop{\longrightarrow}_{n \rightarrow+ \infty}^{\mathcal D}
\biggl(\int_0^t \sigma(X_u,Y_u)
\,d\tilde W_u \sigma(X_u,Y_u),
W_{t_1}, \ldots, W_{t_N} \biggr),
\]
and then apply the results we recalled on stable convergence as we did
in the constant case. To get the previous convergence, one can use the
central limit theorem for triangular arrays of Lindeberg type, stated
for instance in \cite{dacdu}, Theorem~2.8.42.

Another possibility
is to directly use Jacod's stable convergence theorem stated in
Theorem~2.6 of \cite{PV}. Actually in our situation, both theorems require
exactly the same controls (this is not surprising), as soon as one
verifies that the statement of Jacod's theorem extends to a
multi-dimensional setting.

We choose the second solution, and use the
notation in \cite{PV},  Theorem~2.6, so that our $\xi_{n,p}/\sqrt{h_n}$
is equal to their $X_{pn}$.

Conditions (2.6) (martingale increments)
and (2.10) (dependence on $W_.$ only) in \cite{PV} are satisfied.
Condition (2.8) is also satisfied with $v_s=0$ as we already remarked
in the constant case. Here, it amounts to see that
\[
\mathbb{E}_z^0 \bigl[A_{i,j}(n)
\bigl(W^k_{(2p+1)h_n}-W^k_{(2p-1)h_n}\bigr)|
\mathcal F_{(2p-1)h_n} \bigr] = 0
\]
for all triple $(i,j,k)$ where the $A_{i,j}(n)$ are defined in \eqref{offdiag} and \eqref{diag}, which is immediate.

 It thus remains to
check the two conditions
%
%e3.17 #&#
\begin{equation}
\label{eqJac1} \quad\frac{1}{h_n} \sum_{p=1}^{[{t}/({2h_n})]-1}
\mathbb {E}_z^0 \bigl[\bigl(^te_i
\xi_{n,p} e_j\bigr)^2|\mathcal
F_{(2p-1)h_n} \bigr] \mathop{\longrightarrow}_{n \rightarrow+ \infty}^{\mathbb P_z^0} \int
_0^t \theta ^2_{ij}(X_u,Y_u)
\,du,
\end{equation}
for all $i,j=1,\ldots,d$ ($e_l,   l=1,\ldots,d$ being the canonical basis), and
%
%e3.18 #&#
\begin{equation}
\label{eqJac2} \qquad\frac{1}{h_n} \sum_{p=1}^{[{t}/({2h_n})]-1}
\mathbb {E}_z^0 \bigl[|\xi _{n,p}|^2
\BBone_{|\xi_{n,p}|>\varepsilon}|\mathcal F_{(2p-1)h_n} \bigr] \mathop{
\longrightarrow}_{n \rightarrow+ \infty}^{\mathbb
P_z^0} 0 \qquad\mbox{for all } \varepsilon>0,
\end{equation}
where $|\xi|$ denotes the Hilbert--Schmidt norm of the matrix $|\xi|$.

We denote by $u_i=\sigma e_i$, and we use the notation of
Lemma~\ref{lemapprox2}, $U_i(n,s)= \int_{(2p-1)h_n}^{s}
(h_n-|u-2ph_n|)   \,dW^i_u$ and simply $U_i(n)=U_i(n,(2p+1)h_n)$.
Hence,
\[
A_{i,j}(n)=U_i(n) U_j(n) -
\delta_{i,j} \frac{2h_n^3}{3}.
\]
It follows
\begin{eqnarray*}
h_n^4 \mathbb{E}_z^0 \bigl[
\bigl(^te_i \xi_{n,p} e_j
\bigr)^2|\mathcal F_{(2p-1)h_n} \bigr] &=& \mathbb{E}_z^0
\biggl[ \biggl(\sum_{k,l} u_i^k
A_{k,l}(n) u_j^l \biggr)^2\Big|
\mathcal F_{(2p-1)h_n} \biggr]
\\
&=& \mathbb{E} _z^0 \biggl[\sum
_{k,l,k',l'} u_i^k u_j^l
u_i^{k'} u_j^{l'}
A_{k,l}(n) A_{k',l'}(n)|\mathcal F_{(2p-1)h_n} \biggr]
\\
&=& \sum_{k,l,k',l'} u_i^k
u_j^l u_i^{k'}
u_j^{l'} \mathbb {E}_z^0
\bigl[A_{k,l}(n) A_{k',l'}(n)|\mathcal F_{(2p-1)h_n} \bigr].
\end{eqnarray*}
But all conditional expectations are vanishing except those for which
$(k,l)=(k',l')$, in which case it is equal to
\[
(1+\delta_{k,l}) \biggl(\int_{(2p-1)h_n}^{(2p+1)h_n}
\bigl(h_n-|u-2ph_n|\bigr)^2 \,du \biggr)^2
= \frac{4}9 h_n^6 (1+\delta _{k,l}).
\]
Hence,
\begin{eqnarray*}
&&\frac{1}{h_n} \sum_{p=1}^{[{t}/({2h_n})]-1} \mathbb
{E}_z^0 \bigl[\bigl(^te_i
\xi_{n,p} e_j\bigr)^2|\mathcal
F_{(2p-1)h_n} \bigr] \\
&&\qquad= \sum_{k,l=1}^d
(1+\delta_{k,l}) \frac{4h_n}{9} \sum_{p=1}^{[{t}/({2h_n})]-1}
\sigma^2_{i,k}(Z_{(2p-1)h_n}) \sigma^2_{j,l}(Z_{(2p-1)h_n}),
\end{eqnarray*}
and converges to
\[
\sum_{k,l=1}^d (1+\delta_{k,l})
\frac{4}{9} \int_0^t \sigma
^2_{i,k}(Z_u) \sigma^2_{j,l}(Z_u)
\,du.
\]
We get a similar result for $\frac{1}{h_n}   \sum_{p=1}^{[{t}/({2h_n})]-1}   \mathbb{E}_z^0 [(^te_i \xi_{n,p} e_j)(^te_{i'}
\xi_{n,p}
e_{j'})|\mathcal F_{(2p-1)h_n} ]$ for any pairs $(i,j)$,
$(i',j')$. It remains to remark that this increasing process is the one
of
\[
\frac{2}3 \int_0^t
\sigma(Z_u) \,d\tilde{W}_u \sigma(Z_u),
\]
where $\tilde{W}_{.}$ is as in the statement of the proposition.

Finally, \eqref{eqJac2} is immediately checked, using the previous
calculation, Cauchy--Schwarz inequality and Burkholder--Davis--Gundy inequality.
\end{pf}
To conclude the proof of the main theorem, it is enough to apply
Slutsky's theorem since all the error terms converge to $0$ in
probability (recall that Slutsky's theorem also works with stable
convergence (see \cite{PV}, Proposition~2.5).

%s4 #&#
\section{Infinite-horizon estimation}
In the previous section, we dealt with infill estimation. We now
consider that we work with an infinite-horizon design. We aim at
estimating the quantity $\mathbb{E}_\mu(\sigma^2(X_0,Y_0))$, where
$(Z_t:=(X_t,Y_t) \in\mathbb{R}^{2},   t \geq0 )$ is
still governed by (\ref{WuGirsanov}) and $\mu$ is the invariant
measure, supposed to exist. We thus have to introduce some new assumptions:
\begin{itemize}
\item $\mathcal{H}_3$ There exists an (unique) invariant
probability measure $\mu$ and the $\mathbb P_\mu$ stationary process
$Z_.$ is $\alpha$-mixing with rate $\tau$, that is (in our Markovian
situation), there exists a nonincreasing function $\tau$ going to $0$
at infinity such that for all $u\leq s$, all random variables $F,G$
bounded by $1$ s.t. $F$ (resp., $G$) is $\mathcal F_u$ (resp.,
$\mathcal G_s$) measurable where $\mathcal F_u$ (resp., $\mathcal G_s$)
is the $\sigma$-algebra generated by $Z_v$ for $v\leq u$ (resp.,
$v\geq
s$), one has
\[
\operatorname{Cov}_\mu\bigl(F(Z_u) G(Z_s)
\bigr) \leq\tau(u-s).
\]

\item $\mathcal{H}_4$   Define $b(x,y):=- (c(x,y)y+
\nabla V (x)  )$. There exists some $r\geq4$ such that $\mathbb
{E}_\mu
(|b(Z_0)|^r)<+\infty$ and $\int_0^{+\infty}   \tau^{1-(4/r)}(t)\,dt
<+\infty$.
\end{itemize}

These apparently technical assumptions are in a sense ``minimal'' for
applying known results on the central limit theorem for additive
functionals of a diffusion process (see, e.g., \cite{CCG}). We shall
come back later to these assumptions, indicating in the last subsection
of this section, sufficient conditions for them to hold.

We introduce the following estimator:
%
%e4.1 #&#
\begin{equation}
\label{infiniteestimate}
\mathcal{K}_n=\frac{3}2 \frac{1}{(n-1)h_n^3}
\sum_{p=1}^{n-1} \Delta _2 X(p,n)
\otimes\Delta_2 X(p,n),
\end{equation}
where $\Delta_2X(p,n)$ is the\vadjust{\goodbreak} double increment of $X$ defined in
\eqref{acdouble}.

We now state the main result of this section.

%th4.1 #&#
\begin{Theorem}\label{cverg}
Assume that $\mathcal H_0$ up to $\mathcal H_4$ are
satisfied. Assume in addition that
\[
\int_1^{+\infty} t^{-1/2}
\tau^{1/2}(t)\,dt <+\infty.
\]
Let $h_n$ be a sequence going to $0$ such that $nh_n \to
+\infty$ and $nh_n^3 \to0$.

Then, in the stationary regime,
%
%e4.2 #&#
\begin{equation}
\label{regime1} \sqrt{2nh_n} \bigl(\mathcal{K}_n -
\mathbb{E}_{\mu} \sigma^2 (X_0,Y_0)
\bigr)\mathop{\longrightarrow}_{n \rightarrow+ \infty}^{\mathcal{D}} \mathcal{N},
\end{equation}
where $\mathcal N$ is a symmetric random matrix, with centered Gaussian
entries satisfying
\[
\operatorname{Cov}(\mathcal N_{i,j},\mathcal N_{k,l}) =
\frac{1}2 \int_0^{+\infty}
\mathbb{E}_\mu\bigl(\bar{\sigma}^2_{i,j}(Z_0)
\bar{\sigma}^2_{k,l}(Z_s) + \bar {
\sigma}^2_{k,l}(Z_0) \bar{\sigma}^2_{i,j}(Z_s)
\bigr) \,ds,
\]
where $\bar{\sigma}^2(z)=\sigma^2(z)-\mathbb{E}_\mu(\sigma^2(Z_0))$.
\end{Theorem}

%re4.2 #&#
\begin{remark}
In the case where $\sigma$ is constant, this result is useless as the
covariances are all vanishing.
\end{remark}

\begin{pf*}{Proof of Theorem~\protect\ref{cverg}}
From now on, we assume that the assumptions $\mathcal H_0$ up to
$\mathcal H_4$ are satisfied.

Of course since we are looking at the whole time interval up
to infinity, it is no more possible to use Girsanov theory to reduce
the problem to $c=V=0$. Hence, arguing as for the statement of \eqref{ecriture}, and defining $b(x,y):=- (c(x,y)y+ \nabla V (x)
)$, we get
\[
\Delta_2X(p,n)=\int_{(2p-1)h_n}^{(2p+1)h_n}
\bigl(h_n-|u-2ph_n|\bigr) \bigl(\sigma (Z_u)
\,dW_u+b(Z_u)\,du \bigr).
\]
We then define the semimartingale $(H_t,   (2p-1)h_n \leq t
\leq(2p+1)h_n )$ by
\begin{eqnarray*}
dH_t & = & \bigl(h-|t-2ph_n|\bigr)
\sigma(Z_t)\,dW_t+\bigl(h_n-|t-2ph_n|\bigr)b(Z_t)
\,dt,
\\
H_{(2p-1)h_n} & = & 0,
\end{eqnarray*}
so that $\Delta_2X(p,n)=H_{(2p+1)h_n}$. Using It\^{o}'s formula, we then
have
%e4.3 #&#
\begin{eqnarray}
&& \bigl(\Delta_2 X(p,n)\otimes\Delta_2
X(p,n)
\bigr)_{i,j}\nonumber \\
 &&\qquad=\int_{(2p-1)h_n}^{(2p+1)h_n}
\bigl(h_n-|u-2ph_n|\bigr) \bigl(H_u^i
\bigl(\sigma (Z_u) \,dW_u\bigr)^j +
H_u^j \bigl(\sigma(Z_u)
\,dW_u\bigr)^i\bigr)
\nonumber
\\[-8pt]
\label{eqitodelta}
\\[-8pt]
\nonumber
&&\qquad{}+ \int_{(2p-1)h_n}^{(2p+1)h_n} \bigl(h_n-|u-2ph_n|\bigr)
\bigl(H_u^i b^j(Z_u) +
H_u^j b^i(Z_u)\bigr) \,du
\\
&& \qquad{}+ \int_{(2p-1)h_n}^{(2p+1)h_n} \bigl(h_n-|u-2ph_n|\bigr)^2
\sigma^2_{i,j}(Z_u) \,du.
\nonumber
\end{eqnarray}

We have a simple but useful estimate, available for all
$i=1,\ldots, d$, all $k\in\mathbb N$, all $p$ and all $u$ between
$(2p-1)h_n$ and $(2p+1)h_n$
%
%e4.4 #&#
\begin{eqnarray}
&& \mathbb{E}_\mu\bigl(\bigl|H_s^i\bigr|^{2k}
\bigr)
\nonumber
\\
&&\label{eqestimateuse}\qquad\leq C(k) \Vert\sigma\Vert_\infty^{2k} \bigl(s-(2p-1)h_n
\bigr)^k h_n^{2k} \\
&&\qquad\quad{}+ \bigl(s-(2p-1)h_n
\bigr)^{2k} h_n^{2k} \bigl(\mathbb{E}
_\mu \bigl(\bigl|b(Z_0)\bigr|^{2k}\bigr)\bigr).\nonumber
\end{eqnarray}
Indeed, one can first use $(a+b)^{2k}\leq
C(k)(a^{2k}+b^{2k})$, for positive numbers $a,b$ which will be here the
absolute values of the martingale part and of the bounded variation
part.

Then, if $b_u$ is stationary and $h_u$ bounded by $h$,
\begin{eqnarray*}
\mathbb{E}_\mu \biggl( \biggl(\int_0^t
b_u h_u \,du \biggr)^m \biggr) &\leq&
t^{m} h^{m} \mathbb{E}_\mu\bigl(b_0^{m}
\bigr),
\end{eqnarray*}
which can be used with $m=2k$, $t=(s-(2p-1)h_n)$, $b_u=b^i(Z_u)$,
$h_u=(h_n-|u-2ph_n|) \leq2h_n$. This gives the control for the bounded
variation part. Finally, using the Burkholder--Davis--Gundy inequality,
we are reduced to the same control for the martingale part; this time
with $m=k$, $h_u=(h_n-|u-2ph_n|)^2 \leq4h^2_n$ and $|b_u|\leq\Vert \sigma
\Vert_\infty^2$.

Now we can decompose
\[
\mathcal{K}_n - \mathbb{E}_\mu\sigma^2(Z_0)=
\overline{\mathcal {K}}_{n,1}+\overline{\mathcal{K}}_{n,2}
\]
with
\begin{eqnarray*}
\overline{\mathcal{K}}_{n,1}&=&\frac{3}2 \frac{1}{(n-1)h_n^3}
\sum_{p=1}^{n-1} \biggl\{
\Delta_2 X(p,n)\otimes\Delta_2 X(p,n) \\
&&{}- \int
_{(2p-1)h_n}^{(2p+1)h_n}\bigl(h_n-|s-2ph_n|\bigr)^2
\sigma^2(Z_s)\,ds \biggr\}
\end{eqnarray*}
and
\[
\overline{\mathcal{K}}_{n,2}=\frac{3}2 \frac{1}{(n-1)h_n^3}
\sum_{p=1}^{n-1} \int_{(2p-1)h_n}^{(2p+1)h_n}\bigl(h_n-|s-2ph_n|\bigr)^2
\bigl\{ \sigma ^2(Z_s) -\mathbb{E}_\mu
\sigma^2(Z_0) \bigr\}\,ds.
\]
We shall look at both quantities separately, starting with $\overline
{\mathcal{K}}_{n,2}$.
\end{pf*}

%s4.1 #&#
\subsection{Study of $\overline{\mathcal{K}}_{n,2}$}

%le4.3 #&#
\begin{Lemma}\label{lemK2}
There exists some constant $C$ only depending on the bounds of $\sigma$
such that
\[
\mathbb{E}_\mu \bigl\{|\overline{\mathcal{K}}_{n,2}|^2
\bigr\} \leq\frac{C}{n
  h_n} \int_0^{+\infty}
\tau(t) \,dt.
\]
\end{Lemma}

\begin{pf}
\begin{eqnarray*}
&& \frac{4}9 \mathbb{E}_\mu \bigl\{ (\overline
{\mathcal
{K}}_{n,2} )_{i,j}^2 \bigr\} \\
&&\qquad=   \frac{1}{(n-1)^2h_n^6}
\sum_{p,q=1}^{n-1}\int_{(2p-1)h_n}^{(2p+1)h_n}
\int_{(2q-1)h_n}^{(2q+1)h_n}\bigl(h_n-|s-2ph_n|\bigr)^2\\
&&\qquad\quad{}\times\bigl(h_n-|u-2qh_n|\bigr)^2
\mathbb {E}_\mu \bigl\{ \overline{\sigma}^2_{i,j}(Z_s)
\overline{\sigma}^2_{i,j}(Z_u) \bigr\}\,ds\,du
\\
&&\qquad\leq\frac{\Vert\overline{\sigma}^2\Vert^2_{\infty}}{(n-1)^2h_n^6} \sum_{p,q=1}^{n-1}\int
_{(2p-1)h_n}^{(2p+1)h_n}\int_{(2q-1)h_n}^{(2q+1)h_n}\bigl(h_n-|s-2ph_n|\bigr)^2\\
&&\qquad\quad{}\times\bigl(h_n-|u-2qh_n|\bigr)^2
\tau \bigl(|s-u|\bigr) \,du \,ds
\\
&&\qquad\le\frac{C\Vert \overline{\sigma}^2\Vert ^2_\infty}{(n-1)}+\frac
{C\Vert \overline{\sigma}^2\Vert ^2_{\infty
}}{(n-1)^2}\sum_{|p-q|\geq2}
\tau\bigl(2(|p-q|-1)h_n\bigr)\\
&&\qquad \leq\frac
{C\Vert \overline{\sigma}^2
\Vert ^2_\infty}{(n-1)}+
\frac{C\Vert \overline{\sigma}^2\Vert ^2_{\infty
}}{(n-1)}\sum_{k=1}^{n-2} \tau( 2k
h_n)
\\
&&\qquad \leq C\bigl\Vert \overline{\sigma}^2\bigr\Vert ^2_\infty   \biggl(\frac
{1}{n-1} +
\frac
{1}{(n-1)h_n}   \int_0^{+\infty}   \tau(t)   \,dt \biggr)
\end{eqnarray*}
with $C$ some constant. We have used the fact that $\tau$ is
nonincreasing for the final inequality.
\end{pf}

The previous result indicates why the normalization $\sqrt{nh_n}$ has
to be chosen. Now we decompose again
\[
\overline{\mathcal{K}}_{n,2}=\overline{\mathcal {K}}_{n,21}+
\overline {\mathcal{K}}_{n,22}
\]
by decomposing
\[
\sigma^2(Z_s) -\mathbb{E}_\mu
\sigma^2(Z_0)=\sigma^2(Z_s) -
\sigma ^2(Z_{(2p-1)h_n})+\sigma^2(Z_{(2p-1)h_n}) -
\mathbb{E}_\mu\sigma^2(Z_0).
\]
We thus have
\begin{eqnarray*}
\overline{\mathcal{K}}_{n,22} &=& \frac{1}{2(n-1)h_n} \sum
_{p=1}^{n-1} \int_{(2p-1)h_n}^{(2p+1)h_n}
\bigl(\sigma^2(Z_{(2p-1)h_n}) -\mathbb {E}_\mu
\sigma^2(Z_0)\bigr) \,ds
\\
&=& \frac{1}{2(n-1)h_n} \int_{h_n}^{(2n-1)h_n} \bigl(
\sigma^2(Z_s) -\mathbb{E}_\mu
\sigma^2(Z_0)\bigr) \,ds
\\
&& {}+ \frac{1}{2(n-1)h_n} \sum_{p=1}^{n-1}
\int_{(2p-1)h_n}^{(2p+1)h_n} \bigl(\sigma^2(Z_{(2p-1)h_n})
- \sigma^2(Z_s)\bigr) \,ds
\\
&=& \overline{\mathcal{K}}_{n,222}+\overline{\mathcal
{K}}_{n,221}.
\end{eqnarray*}
It follows that
\begin{eqnarray*}
&&\sqrt{2(n-1)h_n}   \overline{\mathcal{K}}_{n,2}\\
&&\qquad=
\frac{1}{\sqrt{2(n-1)h_n}} \int_{h_n}^{(2n-1)h_n} \bigl(
\sigma^2(Z_s) -\mathbb{E}_\mu
\sigma^2(Z_0)\bigr) \,ds \\
&&\qquad\quad{}+ \sqrt{2(n-1)h_n} (
\overline {\mathcal{K}}_{n,221}+\overline{\mathcal{K}}_{n,21}),
\end{eqnarray*}
the first summand being the important term the two others being error
terms. We shall show that these errors terms converge to $0$ in
$\mathbb L^2$. Indeed,
\begin{eqnarray*}
&&\mathbb{E}_\mu \bigl((n-1)^2h^2_n   (\overline
{\mathcal
{K}}_{n,221})_{i,j}^2 \bigr)\\
&&\qquad \leq C \sum_{p,q=1}^{n-1}\int_{(2p-1)h_n}^{(2p+1)h_n}\int_{(2q-1)h_n}^{(2q+1)h_n} \,ds \,du
\\
&&\qquad\quad{}\times\mathbb{E}_\mu \bigl\{\bigl(\sigma^2_{i,j}(Z_s)-
\sigma ^2_{i,j}(Z_{(2p-1)h_n})\bigr) \bigl(
\sigma^2_{i,j}(Z_u)-\sigma ^2_{i,j}(Z_{(2q-1)h_n})
\bigr) \bigr\},
\end{eqnarray*}
so that, as for the proof of Lemma \ref{lemK2}, what has to be done is
to control
\[
\operatorname{Cov}\bigl(\sigma^2_{i,j}(Z_s)-
\sigma ^2_{i,j}(Z_{(2p-1)h_n}),\sigma
^2_{i,j}(Z_u)-\sigma^2_{i,j}(Z_{(2q-1)h_n})
\bigr).
\]

The problem is that, if we use the $\alpha$-mixing we will not improve
upon the bound in the previous lemma, since the uniform bound of these
variables is still of order a constant. However, for Markov diffusion
processes, one can show (see, e.g., \cite{CCG},  Lemma~4.2 and
Lemma~5.1, or \cite{DOU}, Chapter~1, but the latter result also easily
follows from the Riesz--Thorin interpolation theorem) the following.

%le4.4 #&#
\begin{Lemma}\label{mixfaible}
Let $F$ and $G$ be as in the definition of the $\alpha$-mixing except
that they are not bounded. Assume that $F \in\mathbb L^r(\mathbb
{E}_\mu)$ and
$G \in\mathbb L^a(\mathbb{E}_\mu)$ for some $r$ and $a$ larger than
or equal
to $2$. Then
\begin{eqnarray*}
&&\operatorname{Cov}_\mu(F,G) \\
&&\qquad\leq C \min \bigl(\tau^{({r-2})/({2r})}(s-u)
\Vert F\Vert_{\mathbb L^r} \Vert G\Vert_{\mathbb L^2} ; \tau^{({a-2})/({2a})}(s-u)
\Vert F\Vert_{\mathbb L^2} \Vert G\Vert_{\mathbb L^a} \bigr),
\end{eqnarray*}
for some constant $C$ depending on $a$ and $r$ only. One also has
\[
\operatorname{Cov}_\mu(F,G) \leq C \tau^{({r-2})/({2r})}\bigl((s-u)/2
\bigr) \tau ^{({a-2})/({2a})}\bigl((s-u)/2\bigr) \Vert F\Vert_{\mathbb L^r}
\Vert G\Vert_{\mathbb L^a},
\]
for some constant $C$ depending on $a$ and $r$ only.
\end{Lemma}

Choosing\vspace*{1.5pt} $F=\sigma^2_{i,j}(Z_s)-\sigma^2_{i,j}(Z_{(2p-1)h_n})$ and
$G=\sigma^2_{i,j}(Z_u)-\sigma^2_{i,j}(Z_{(2q-1)h_n})$, $r=a$, we see
that what we have to do is to get a nice upper bound for $\mathbb
{E}_\mu
(|F|^r)$. But
\[
\bigl|\sigma^2_{i,j}(Z_s)-\sigma^2_{i,j}(Z_{(2p-1)h_n})\bigr|
\leq K |Z_s -Z_{(2p-1)h_n}|,
\]
where $K$ only depends on $\sigma$ and its first derivatives. Using
Burkholder--Davis--Gundy inequality, we thus have
\begin{eqnarray*}
\mathbb{E}_\mu\bigl(|F|^r\bigr)&\leq& C
\bigl(h_n^{r/2} + h_n^r
\mathbb{E}_\mu \bigl(\bigl|b(Z_0)\bigr|^r\bigr)\bigr).
\end{eqnarray*}
It follows that, provided $\mathbb{E}_\mu(|b(Z_0)|^r)<+\infty$,
\begin{eqnarray*}
&&\operatorname{Cov}\bigl(\sigma^2_{i,j}(Z_s)-
\sigma ^2_{i,j}(Z_{(2p-1)h_n}),\sigma
^2_{i,j}(Z_u)-\sigma^2_{i,j}(Z_{(2q-1)h_n})
\bigr)\\
&&\qquad\leq C h_n \tau ^{1-(2/r)}\bigl(\bigl(|p-q|-1\bigr)h_n
\bigr),
\end{eqnarray*}
so that finally, as in the proof of Lemma~\ref{lemK2} we get
%
%e4.5 #&#
\begin{equation}
\label{eqerrorK2}
\mathbb{E}_\mu \bigl((n-1)h_n (\overline{
\mathcal {K}}_{n,221})_{i,j}^2 \bigr)\leq C
h_n \biggl(1+\int_0^{+\infty}
\tau^{1-(2/r)}(t) \,dt \biggr).
\end{equation}
Exactly in the same way, we obtain the same result replacing $\overline
{\mathcal{K}}_{n,221}$ by $\overline{\mathcal{K}}_{n,21}$.

It remains
to look at
\[
\frac{1}{\sqrt{2(n-1)h_n}} \int_{h_n}^{(2n-1)h_n} \bigl(
\sigma^2(Z_s) -\mathbb{E}_\mu
\sigma^2(Z_0)\bigr) \,ds.
\]
The asymptotic behavior of such additive functionals of stationary
Markov processes has been extensively studied. For simplicity, we refer
to the recent \cite{CCG} for an overview and a detailed bibliography.
In particular, Section~4 of this reference contains the following
result (essentially due to Maxwell and Woodroofe), provided $\int_1^{+\infty} t^{-1/2}   \tau^{1/2}(t)   \,dt < +\infty$, the previous
quantity converges in distribution to a centered Gaussian random
variable, as soon as $nh_n$ goes to infinity. The calculation of the
covariance matrix of these variables is done as in \cite{CCG}. We have
thus obtained the following.

%pr4.5 #&#
\begin{Proposition}\label{propK2}
Assume\vspace*{1pt} that $\mathcal H_0$ up to $\mathcal H_4$ are satisfied. Assume
in addition that $\int_1^{+\infty}   t^{-1/2}   \tau^{1/2}(t)\,dt
<+\infty$. Let $h_n$ be a sequence going to $0$ such that $nh_n \to
+\infty$.

Then, in the stationary regime,
$\sqrt{2(n-1)h_n}   \overline{\mathcal{K}}_{n,2}$ converges in
distribution to a symmetric random matrix $\mathcal N$, with centered
Gaussian entries satisfying
\[
\operatorname{Cov}(\mathcal N_{i,j},\mathcal N_{k,l}) =
\frac{1}2 \int_0^{+\infty}
\mathbb{E}_\mu\bigl(\bar{\sigma}^2_{i,j}(Z_0)
\bar{\sigma}^2_{k,l}(Z_s) + \bar {
\sigma}^2_{k,l}(Z_0) \bar{\sigma}^2_{i,j}(Z_s)
\bigr) \,ds.
\]
\end{Proposition}

%s4.2 #&#
\subsection{Study of $\overline{\mathcal{K}}_{n,1}$}

%le4.6 #&#
\begin{Lemma}\label{lemK1}
Assume that for some $k\in\mathbb N^*$, $\mathbb{E}_\mu
(|b(Z_0)|^{4k})<+\infty
$ and that $\int_0^{+\infty}   \tau^{1-(1/k)}(t)\,dt <+\infty$.

Then
there exists some constant $C(k)$ such that for all $i,j=1,\ldots,d$,
\[
\operatorname{Var}_\mu \bigl[ (\overline{\mathcal
{K}}_{n,1} )_{i,j} \bigr] \leq\frac{C(k)}{n}.
\]
Hence,
\[
\operatorname{Var}_\mu \bigl[\sqrt{nh_n} (\overline {
\mathcal {K}}_{n,1} )_{i,j} \bigr] \to 0.
\]
\end{Lemma}
\begin{pf}
We write
\[
\Delta_2 X(p,n)\otimes\Delta_2 X(p,n) - \int
_{(2p-1)h_n}^{(2p+1)h_n}\bigl(h_n-|s-2ph_n|\bigr)^2
\sigma^2(Z_s)\,ds =M_{p,n}+V_{p,n},
\]
where\vspace*{1pt} $M_.$ (resp., $V_.$) denotes the martingale (resp., bounded
variation) part. As usual, we use $\bar V$ for the centered $V -\mathbb
{E}_\mu
(V)$. Hence,
\begin{eqnarray*}
&&\frac{4}9 (n-1)^2 h_n^6
\operatorname{Var}_\mu \bigl[ (\overline {\mathcal
{K}}_{n,1} )_{i,j} \bigr] \\
&&\qquad= \sum_{p,q=1}^{n-1}
\mathbb {E}_\mu \bigl(M^{i,j}_{p,n}M^{i,j}_{q,n}+M^{i,j}_{p,n}
\bar V^{i,j}_{q,n}+\bar V^{i,j}_{p,n}M^{i,j}_{q,n}+
\bar V^{i,j}_{p,n}\bar V^{i,j}_{q,n}\bigr).
\end{eqnarray*}
A lot of terms of this sum are vanishing, so that we get
\begin{eqnarray*}
\frac{4}9 (n-1)^2 h_n^6
\operatorname{Var}_\mu \bigl[ (\overline {\mathcal
{K}}_{n,1} )_{i,j} \bigr] &=& \sum
_{p=1}^{n-1} \mathbb{E}_\mu \bigl(
\bigl(M^{i,j}_{p,n}\bigr)^2 + 2 \bar
V^{i,j}_{p,n}M^{i,j}_{p,n}+ \bigl(\bar
V^{i,j}_{p,n}\bigr)^2\bigr)
\\
&& {}+ \sum_{p> q=1}^{n-1} \mathbb{E}_\mu
\bigl(\bar V^{i,j}_{p,n}M^{i,j}_{q,n}+ 2
\bar V^{i,j}_{p,n}\bar V^{i,j}_{q,n}\bigr).
\end{eqnarray*}
Using stationarity and \eqref{eqestimateuse}, we get
\begin{eqnarray*}
\mathbb{E}_\mu\bigl(\bigl(M^{i,j}_{p,n}
\bigr)^2\bigr)&=&\int_0^{2h_n}
\bigl(h_n-|u-h_n|\bigr)^2\\
&&{}\times \mathbb{E} _\mu
\bigl(\sigma_{i,i}^2(Z_s) \bigl(H^j_s
\bigr)^2+\sigma_{j,j}^2(Z_s)
\bigl(H^i_s\bigr)^2+2 \sigma
_{i,j}^2(Z_s) H^j_sH^i_s
\bigr) \,ds
\\
&\leq& C h_n^6 \bigl(1+h_n \mathbb{E}
_\mu \bigl(\bigl|b(Z_0)\bigr|^2\bigr)\bigr).
\end{eqnarray*}
Similarly,
\begin{eqnarray*}
\mathbb{E}_\mu\bigl(\bigl(V^{i,j}_{p,n}
\bigr)^2\bigr)&=&\mathbb{E}_\mu \biggl[ \biggl(\int
_{0}^{2h_n} \bigl(h_n-|u-h_n|\bigr)
\bigl(H_u^i b^j(Z_u) +
H_u^j b^i(Z_u)\bigr) \,du
\biggr)^2 \biggr]
\\
&\leq& C h_n^3 \int_0^{2h_n}
\mathbb{E}_\mu\bigl(\bigl|b(Z_u)\bigr|^2
\bigl(H_u^i\bigr)^2\bigr) \,du
\\
&\leq& C h_n^7 \bigl(\mathbb{E}_\mu
\bigl(\bigl|b(Z_0)\bigr|^4\bigr)\bigr)^{1/2}
\bigl(1+h_n \bigl(\mathbb{E} _\mu \bigl(\bigl|b(Z_0)\bigr|^4
\bigr)\bigr)^{1/2}\bigr).
\end{eqnarray*}
It follows that
\[
\sum_{p=1}^{n-1} \mathbb{E}_\mu
\bigl(\bigl(M^{i,j}_{p,n}\bigr)^2 + 2 \bar
V^{i,j}_{p,n}M^{i,j}_{p,n}+ \bigl(\bar
V^{i,j}_{p,n}\bigr)^2\bigr) \leq C (n-1)
h_n^6.
\]
Exactly in the same way one obtains that, for $k\in\mathbb
N^*$, provided $\mathbb{E}_\mu(|b(Z_0)|^{2k})<+\infty$,
\[
\mathbb{E}_\mu\bigl(\bigl|M^{i,j}_{p,n}\bigr|^{2k}
\bigr)\leq C(k) h_n^{6k}
\]
and provided $\mathbb{E}_\mu(|b(Z_0)|^{4k})<+\infty$,
\[
\mathbb{E}_\mu\bigl(\bigl|V^{i,j}_{p,n}\bigr|^{2k}
\bigr)\leq C(k) h_n^{7k}.
\]
Again we shall use Lemma~\ref{mixfaible} to control
\[
\mathbb{E}_\mu\bigl(\bar V^{i,j}_{p,n}M^{i,j}_{q,n}
\bigr)=\operatorname {Cov}_\mu \bigl(V^{i,j}_{p,n},M^{i,j}_{q,n}
\bigr) \quad\mbox{and}\quad \mathbb{E}_\mu \bigl(\bar V^{i,j}_{p,n}
\bar V^{i,j}_{q,n}\bigr)=\operatorname{Cov}_\mu
\bigl(V^{i,j}_{p,n},V^{i,j}_{q,n}\bigr),
\]
and we obtain
\[
\operatorname{Cov}_\mu\bigl(V^{i,j}_{p,n},M^{i,j}_{q,n}
\bigr) \leq C h_n^6 \tau ^{(k-1)/k}\bigl((p-q-1)/2
\bigr)
\]
and
\[
\operatorname{Cov}_\mu\bigl(V^{i,j}_{p,n},V^{i,j}_{q,n}
\bigr) \leq C h_n^6 \tau ^{(k-1)/k}\bigl((p-q-1)/2
\bigr)
\]
provided, respectively, $\mathbb{E}_\mu(|b(Z_0)|^{2k})<+\infty$ and
$\mathbb{E}_\mu
(|b(Z_0)|^{4k})<+\infty$.

We have thus obtained
\[
\sum_{p> q=1}^{n-1} \mathbb{E}_\mu
\bigl(\bar V^{i,j}_{p,n}M^{i,j}_{q,n}+ 2
\bar V^{i,j}_{p,n}\bar V^{i,j}_{q,n}\bigr)
\leq C (n-1) h_n^6 \int_0^{+\infty}
\tau^{(k-1)/k}(t) \,dt,
\]
so that gathering all previous estimates we get the result.
\end{pf}

It remains to bound the expectation of $ (\overline
{\mathcal{K}}_{n,1} )_{i,j}$. But
\begin{eqnarray*}
&&  \mathbb{E}_\mu \bigl[ (\overline{\mathcal
{K}}_{n,1}
)_{i,j} \bigr]
\\
&&\qquad= \frac{3}{2(n-1)h_n^3}\\
&&\qquad\quad{}\times \sum_{p=1}^{n-1}
\mathbb{E} _\mu \biggl[\int_{(2p-1)h_n}^{(2p+1)h_n}
\bigl(h_n-|u-2ph_n|\bigr) \bigl(H_u^i
b^j(Z_u)+H_u^j
b^i(Z_u)\bigr) \,du \biggr]
\\
&&\qquad= \frac{3}{2(n-1)h_n^3} (A_{n,1}+A_{n,2})
\end{eqnarray*}
with
\begin{eqnarray*}
A_{n,1}&=&\sum_{p=1}^{n-1}
\mathbb{E}_\mu \biggl[\int_{(2p-1)h_n}^{(2p+1)h_n}
\bigl(h_n-|u-2ph_n|\bigr) \bigl(H_u^i
b^j(Z_{(2p-1)h_n})\\
&&{}+H_u^j
b^i(Z_{(2p-1)h_n})\bigr) \,du \biggr],
\end{eqnarray*}
and
\begin{eqnarray*}
A_{n,2} &=& \sum_{p=1}^{n-1}   \mathbb{E}_\mu
\Biggl[\int_{(2p-1)h_n}^{(2p+1)h_n}   \bigl(h_n-|u-2ph_n|\bigr)
\bigl(H_u^i \bigl(b^j(Z_u)-b^j(Z_{(2p-1)h_n})
\bigr)\\
&&{}+H_u^j \bigl(b^i(Z_u)-b^i(Z_{(2p-1)h_n})
\bigr)\bigr) \,du \Biggr].
\end{eqnarray*}
$A_{n,2}$ can be studied exactly as we did before because
$b^j(Z_u)-b^j(Z_{(2p-1)h_n})$ is centered. To be more precise, instead
of calculating $A_{n,2}$ we look at the $\mathbb L^2$ norm of the
random variable
\begin{eqnarray*}
&&\sum_{p=1}^{n-1} \int_{(2p-1)h_n}^{(2p+1)h_n}
\bigl(h_n-|u-2ph_n|\bigr) \bigl(H_u^i
\bigl(b^j(Z_u)-b^j(Z_{(2p-1)h_n})
\bigr)\\
&&\qquad{}+H_u^j \bigl(b^i(Z_u)-b^i(Z_{(2p-1)h_n})
\bigr)\bigr) \,du
\end{eqnarray*}
which is, thanks to the centering property, similar to the quantities
we have studied in the proof of Lemma~\ref{lemK1}, that is, we can use
the mixing property for the \emph{covariances}. It follows that $\sqrt
{nh_n}   \frac{A_{n,2}}{nh_n^3}$ goes to $0$.

Finally, using the semimartingale decomposition of $H_u$,
\begin{eqnarray*}
A_{n,1} &=& \sum_{p=1}^{n-1} \int
_{(2p-1)h_n}^{(2p+1)h_n} \bigl(h_n-|u-2ph_n|\bigr)
\mathbb{E}_\mu\bigl(H_u^i
b^j(Z_{(2p-1)h_n})\\
&&{}+H_u^j
b^i(Z_{(2p-1)h_n})\bigr) \,du
\\
&=& \sum_{p=1}^{n-1} \int
_{(2p-1)h_n}^{(2p+1)h_n} \int_{(2p-1)h_n}^{u}
\bigl(h_n-|u-2ph_n|\bigr)\bigl(h_n-|v-2ph_n|\bigr)
\\
&&{}\times\mathbb{E}_\mu\bigl(b^i(Z_v)
b^j(Z_{(2p-1)h_n})+b^j(Z_v)
b^i(Z_{(2p-1)h_n})\bigr) \,dv\,du
\end{eqnarray*}
so that
\[
|A_{n,1}|\leq C nh_n^4 \bigl(
\mathbb{E}_\mu\bigl(\bigl|b(Z_0)\bigr|^2\bigr)
\bigr)^2.
\]
Hence, $\sqrt{nh_n}   \frac{A_{n,1}}{nh_n^3}$ goes to $0$, provided
$nh_n^3 \to0$. This completes the proof of the theorem.
%\end{pf}

%s4.3 #&#
\subsection{The \texorpdfstring{$\sigma$}{sigma} constant case}

As we already remarked, if $\sigma(x,y)$ is constant, $\overline
{\mathcal{K}}_{n,2} = 0$. The good normalization is then $\sqrt n$.
Indeed, in the previous proof we did not use the full strength of the
bound
\[
\mathbb{E}_\mu\bigl(\bigl|V^{i,j}_{p,n}\bigr|^{2k}
\bigr)\leq C(k) h_n^{7k},
\]
furnishing some $h_n^{7/2}$ instead of a $h_n^3$ each time a bounded
variation term appears. Hence, all terms will go to $0$ except the two
remaining terms:
\begin{itemize}
\item $\sqrt{n}   \frac{A_{n,1}}{nh_n^3} \leq C   (\mathbb{E}_\mu
(|b(Z_0)|^2))^2   \sqrt n   h_n$ for which we need $n h_n^2 \to0$,
\item and the remaining martingale term
\begin{eqnarray*}
&&\int_{(2p-1)h_n}^{(2p+1)h_n} \bigl(h_n-|u-2ph_n|\bigr)
\\
&&\qquad{}\times\int_{(2p-1)h_n}^{u} \bigl(h_n-|s-2ph_n|\bigr)
\bigl(\sigma(Z_s) \,dW_s\bigr)^i \bigl(
\sigma(Z_u) \,dW_u\bigr)^j
\end{eqnarray*}
in \eqref{eqitodelta}.
\end{itemize}
But since $\sigma$ is constant, this is exactly the martingale term we
encountered in Section~\ref{subsecons}. We thus have obtained the following.

%th4.7 #&#
\begin{Theorem}\label{cvergconstant}
Assume that $\mathcal H_0$ up to $\mathcal H_4$ are satisfied
and that $\sigma$ is constant.

Let $h_n$ be a sequence going to $0$ such that $nh_n \to
+\infty$ and $nh_n^2 \to0$.

Then, in the stationary regime,
%
%e4.6 #&#
\begin{equation}
\label{regime11} \sqrt{n} \bigl(\mathcal{K}_n - \sigma^2
\bigr) \mathop{\longrightarrow}_{n\rightarrow+ \infty}^{\mathcal{D}} \sigma
\mathcal{N}_{(d,d)} \sigma,
\end{equation}
where $\mathcal{N}_{(d,d)}$ is as in Lemma~\ref{cvas}.
\end{Theorem}

%s4.4 #&#
\subsection{About $\mathcal{H}_3$ and $\mathcal{H}_4$}
As we promised, we come back to the conditions $\mathcal{H}_3$ and
$\mathcal{H}_4$. Actually, in full generality, very few are known. All
known results amount to the existence of some Lyapunov function (see,
e.g., \cite{Wu}, Theorem~2.4), that is, some nonnegative function
$\psi$ satisfying  $-L\psi\geq\lambda  \psi$ at infinity for some
$\lambda
>0$. In this case, $\tau$ has an exponential decay and the invariant
measure exponential moments, so that $\mathcal H_3$ and $\mathcal H_4$
are satisfied provided $b$ has some polynomial growth. General (and not
really tractable) conditions for the existence of $\psi$ are discussed
in \cite{Wu},  Sections~3 and 4. One can also relax the Lyapunov control
as in \cite{DFG}.

Tractable conditions are only known when $\sigma$ is constant. They are
recalled in \cite{CLP} (see hypotheses $\mathcal H_1$ and $\mathcal
H_2$ therein, based on \cite{Wu} and \cite{BCG}). Mainly, one has to
assume that $c$ and $V$ have at most polynomial growth and that
$<x,\nabla V(x)>$ is positive enough at infinity, for instance,
\[
\bigl\langle x,\nabla V(x) \bigr\rangle \geq\lambda|x|
\]
at infinity.

%s5 #&#
\section{Fluctuation-dissipation relation and Langevin dynamics}\label
{Langevin}
In this section, we focus on Langevin equations, satisfying the
so-called fluct\-uation-dissipation relation. The motivation for the
study of such dynamics comes from modelling interaction of a subsystem
with its environment. Several derivations can be found in the
literature, among others the Ehrenfest dynamics (see, e.g., \cite
{BoBo,SzSz,HoSc}) or the nonlinear Kac--Zwanzig heat bath models
(see, e.g., \cite{StuartKZ}).

We now propose in this section an estimation procedure for Langevin
dynamics satisfying the so-called fluctuation-dissipation relation,
that is, dynamics described by equation (\ref{LangevinEquations}) in
the \hyperref[sec1]{Introduction}. Following Remark~\ref{remsig}, (\ref
{LangevinEquations}) can be written as
\[
\cases{
dX_t = Y_t \,dt,
\vspace*{3pt}\cr
dY_t  = \sqrt{2 \beta^{-1}} s(X_t)
\,dW_t-\bigl(s(X_t)s^*(X_t)Y_t+
\nabla V(X_t)\bigr) \,dt.}
\]
The study of such systems is of great interest for understanding
molecular dynamics (see references previously cited, as well as
\cite{HoSc} and the references therein).

As already mentioned in the \hyperref[sec1]{Introduction}, it is proved under
assumptions $\mathcal{H}_0$ and $\mathcal{H}_1$ of Section~\ref{sectools}, that the solution of (\ref{LangevinEquations}) is
exponentially ergodic with invariant probability measure proportional\vspace*{2pt}
to the Boltzmann distribution $\exp(-\beta H(x,y))$, where
$H(x,y)=\frac{1}2|y|^2+V(x)$ and $\beta$ is inversely proportional to the
temperature (see, e.g., \cite{Matti1,Tony}). In what follows, we denote
by $p_s^\beta(x,y)$ the density of the invariant measure. We shall now
propose an estimation procedure for the parameters associated to the
system described by (\ref{LangevinEquations}).

First, we consider the estimation of the diffusion term. Under
assumptions $\mathcal{H}_0$ and $\mathcal{H}_1$, the results of the
Theorem~\ref{thmmainvarinfill} still hold ($\mathcal{H}_0$ indeed
implies $\mathcal{H}_2$ if the fluctuation-dissipation relation is
satisfied). We thus get
\[
\mathcal{QV}_{h_n}=\frac{1}{h_n^2}\sum
_{p=1}^{[{t}/({2h_n})]-1}\Delta_2 X(p,n)\otimes
\Delta_2 X(p,n)\mathop{\longrightarrow}_{n\rightarrow+ \infty}^{\mathbb
{P}_z}
\frac{2}{3\beta}\int_0^t s(X_s)s^*(X_s)
\,ds
\]
and
\[
\sqrt{\frac{1}{h_n}} \biggl(\mathcal{QV}_{h_n}(t)-
\frac{2}{3\beta} \int_0^t s
(X_s)s^*(X_s)\,ds \biggr) \mathop{\longrightarrow}_{n \rightarrow+ \infty}^{\mathcal S}
\frac{4}{3\beta} \int_0^t s(X_s)
\,d\tilde W_s s(X_s).
\]

The infinite horizon setting can also be considered. Ergodic properties
in \cite{Matti1} and in \cite{Tony} for nonperiodic potentials imply
assumption $\mathcal{H}_3$ with an exponential rate. Moreover,\vspace*{1pt} defining
$g(x,y)=-(s(x)s^*(x)y+\nabla V(x))$, we get from the form of the
invariant density $p_s^\beta(x,y)$: $\mathbb{E}_{\mu
}(|g(Z_0)|^r)<\infty$ for
all $r>0$. Besides, we also have $\int_1^\infty t^{-1/2}\tau
^{1/2}(t)\,dt<\infty$. Thus, assumption $\mathcal{H}_4$ is satisfied and
the results of Theorem~\ref{cverg} still holds:

Let $h_n$ be a sequence going to $0$ such that $nh_n \to
+\infty$ and $nh_n^3 \to0$.

Then, recalling that
\[
\mathcal K_n=\frac{3}2\frac{1}{(n-1)h_n^3}\sum
_{p=1}^{n-1}\Delta _2X(p,n)\otimes
\Delta_2X(p,n),
\]
we have
\begin{eqnarray*}
&&\sqrt{2nh_n} \biggl(\mathcal{K}_n -
\frac{2}\beta\int s(x)s^*(x)p_s^\beta (x,y)\,dx\,dy
\biggr) \mathop{\longrightarrow}_{n \rightarrow+ \infty}^{\mathcal{D}} \mathcal{N},
\end{eqnarray*}
where $\mathcal N$ is a symmetric random matrix, with centered Gaussian
entries satisfying
\[
\operatorname{Cov}(\mathcal N_{i,j},\mathcal N_{k,l}) =
\frac{1}2 \int_0^{+\infty}
\mathbb{E}_\mu\bigl(\bar{s}^2_{i,j}(X_0)
\bar{s}^2_{k,l}(X_s) + \bar
{s}^2_{k,l}(X_0) \bar{s}^2_{i,j}(X_s)
\bigr) \,ds,
\]
where $\bar{s}^2(x)=\frac{2}\beta(s(x)s^*(x)-\mathbb{E}_\mu
(s(X_0)s^*(X_0)))$.
Notice that
\[
\frac{2}\beta\mathbb{E}_\mu\bigl(s(X_0)s^*(X_0)
\bigr)=\frac{2}\beta\int s(x)s^*(x)p_s^\beta(x,y)
\,dx\,dy.
\]
Although we have been able to estimate the quadratic variation in both
cases there are parameters that remain undetermined.

This leads us to consider a more general estimation taking into account
our two \cite{CLP2} and \cite{CLP} previous articles. However, to get
easier computations, we shall only consider here the case when the two
coordinates of the process are observed. In this case, the computations
below are simple adaptations of what we have done in our previous
works. The extension to partial observations is not as immediate, and
requires to re-write a large part of these works, but following closely
the same lines of reasoning. This job cannot be done here.

We know that $p_s^\beta(x,y)=\mathbf C(\beta)e^{-\beta(|y|^2+V(x))}$
then $\frac{\nabla_x p_s^\beta(x,y)}{p_s^\beta(x,y)}=-\beta\nabla
V(x)$. To estimate this last quantity, let $K$ be a convolution kernel
with bounded support satisfying $\int K(x,y)\,dx\,dy=1$ and verifying that
there exists an integer $m>0$ such that for all nonconstant polynomial
$P(x,y)$ of degree less than or equal to $m$, $\int P(x,y)
K(x,y)\,dx\,dy=0$. Let $h_n$, $b_{1,n}$ and $b_{2,n}$ be sequences
satisfying the hypothesis (i), (ii), (iii) and (iv) of
Theorem~3.3 in \cite{CLP}. Then we introduce the following estimators:
\begin{eqnarray*}
\tilde p_s(x,y)&=&\frac{1}{nb^d_{1,n}b^d_{2n}}\sum_{i=1}^n
K\biggl(\frac
{x-X_{ih_n}}{b_{1,n}},\frac{y-Y_{ih_n}}{b_{2,n}}\biggr),
\\
\nabla_x\tilde p_s(x,y)&=&\frac{1}{nb^{d+1}_{1,n}b^d_{2n}}\sum
_{i=1}^n \nabla_xK\biggl(
\frac{x-X_{ih_n}}{b_{1,n}},\frac{y-Y_{ih_n}}{b_{2,n}}\biggr).
\end{eqnarray*}
The candidate for estimating $-\beta\nabla V(x)$ will be
$\frac{\nabla_x\tilde p_s(x,y)}{\tilde p_s(x,y)}$ whose consistence in
probability and asymptotic normality is derived below. Let us write
\begin{eqnarray*}
 \mathcal A_n(x,y) &:=& \biggl(\frac{\nabla_x\tilde
p_s(x,y)}{\tilde p_s(x,y)}+\beta\nabla V(x)\biggr)
\\
&=&\biggl(\frac{\nabla_x\tilde p_s(x,y)}{\tilde p_s(x,y)}-\frac{\nabla_x
p_s^\beta(x,y)}{\tilde p_s(x,y)}\biggr)+\biggl(
\frac{\nabla_x p_s^\beta
(x,y)}{\tilde p_s(x,y)}+\beta\nabla V(x)\biggr)
\\
&=&\frac{1}{\tilde p_s(x,y)}\bigl(\nabla_x\tilde p_s(x,y)-
\nabla_x p_s^\beta (x,y)\bigr)\\
&&{}-
\frac{\nabla_x p_s^\beta(x,y)}{\tilde p_s(x,y)p_s^\beta
(x,y)}\bigl(\tilde p_s(x,y)-p_s^\beta(x,y)
\bigr).
\end{eqnarray*}
Following the proof of Theorem~3.3 in \cite{CLP}, we can prove that the
second term in the last equality above is $O_{\mathbb{P}_z}(\sqrt{n
b^d_{1,n}b^d_{1,n}})$. Recalling\vspace*{1pt} that $\mathcal D$ denotes the
convergence in distributions of probability measures, we have
\begin{eqnarray*}
&& \mathcal D\lim_{n\to\infty}\sqrt{n
b^{(d+2)}_{1,n}b^d_{2,n}}\mathcal A_n(x,y)
\\
&&\qquad=\mathcal D\lim_{n\to\infty}\sqrt{n b^{(d+2)}_{1,n}b^d_{2,n}}
\frac{1}{\tilde p_s(x,y)}\bigl(\nabla_x\tilde p_s(x,y)-
\nabla_x p_s^\beta(x,y)\bigr)
\\
&&\qquad=\frac{1}{p_s^\beta(x,y)}\mathcal D\lim_{n\to\infty}\sqrt{n
b^{(d+2)}_{1,n}b^d_{2,n}}\bigl(
\nabla_x\tilde p_s(x,y)-\nabla_x
p_s^\beta(x,y)\bigr).
\end{eqnarray*}
The last equality above is a consequence of Slutsky's theorem.

We now
sketch the proof of the convergence in distribution of
$\mathcal R_n:=\sqrt{n b^{(d+2)}_{1,n}b^d_{2,n}}(\nabla_x\tilde
p_s(x,y)-\nabla_x p_s^\beta(x,y))$. Let us denote by $\partial x_l$ the
partial derivative with respect to the $l$th coordinate.
Using as a tool the computation of covariances for sums of $\alpha
$-mixing random variables, we get
\[
\operatorname{Cov}\bigl(\partial_{x_j}\tilde p(x,y),
\partial_{x_l}\tilde p(x,y)\bigr)=O\biggl(\frac
{p_s^\beta(x,y)}{n b^{(d+2)}_{1,n}b^d_{2,n}}\int\partial
_{x_j}K(u,v)\partial_{x_l}K(u,v)\,du\,dv\biggr).
\]
Hence, and as a consequence that the kernel $K$ has bounded support we have
\[
\bigl(n b^{(d+2)}_{1,n}b^d_{2,n}\bigr)
\operatorname{Cov}\bigl(\partial_{x_j}\tilde p(x,y),\partial
_{x_l}\tilde p(x,y)\bigr)\to\delta_{jl}p_s^\beta(x,y)
\int\bigl(\partial _{x_j}K(u,v)\bigr)^2\,du\,dv,
\]
where the $\delta_{ij}$'s stand for the Kronecker symbols. The random
sequence $\mathcal R_n$ is a sum of a triangular array of $\alpha$-mixing random vectors of $\mathbb{R}^d$. It is straightforward to extend
the results of Theorem~4.3 in \cite{CLP}, via the Cram\'er--Wald device,
to random vectors. Thus, defining $\mathbf D(x,y)=(d_{ij}(x,y))$ as a
diagonal matrix, and if the sequences $h_n$, $b_{1,n}$ and $b_{2,n}$
satisfy the hypothesis of Theorem~3.3 in \cite{CLP}, we get
\[
\mathcal R_n\stackrel{\mathcal D}\to\mathcal N\bigl(0,\mathbf D(x,y)
\bigr)\qquad \mbox{where } d_{jj}(x,y)=p_s^\beta(x,y)
\int\bigl(\partial_{x_j}K(u,v)\bigr)^2\,du\,dv.
\]
We finally obtain
\[
\sqrt{n b^{(d+2)}_{1,n}b^d_{2,n}}
\mathcal A_n(x,y)\stackrel{\mathcal D}\to N\bigl(0,\mathbf L(x,y)\bigr)
\]
where $\mathbf{L}$ is also diagonal with
$l_{jj}(x,y)=\frac{1}{p_s^\beta(x,y)}\int(\partial_{x_j}K(u,v))^2\,du\,dv$.

Now we consider the drift estimation. We can estimate the function
$g(x,y)=-[s(x)s^*(x)y+\nabla V(x)]$. We recall that we only consider
the case of complete observations and for simplicity we consider now
$d=1$. We define the Naradaya--Watson estimator
\[
H_n(x,y)=\frac{1}{(n-1)b^d_{1,n}b^d_{2,n}}\sum_{i=1}^{n-1}K
\biggl(\frac
{x-X_{ih_n}}{b_{1,n}},\frac{y-Y_{ih}}{b_{2,n}}\biggr)\frac
{Y_{(i+1)h_n}-Y_{ih_n}}{h_n}.
\]
The following approximation
%
%e5.1 #&#
\begin{eqnarray}
\label{approx}
&& Y_{(i+1)h_n}-Y_{ih_n}\approx s (X_{ih_n})
(W_{(i+1)h_n}-W_{ih_n})+g(X_{ih_n},Y_{ih_n})h_n,
\end{eqnarray}
permits to obtain
\[
\mathbb{E}\bigl[H_n(x,y)\bigr]\to g(x,y)p^\beta_s(x,y).
\]
We recall that we provide here only a flavor of the proof, and not a
rigorous justification for each point. We now write
%$((n-1)b^d_{1,n}b^d_{2,n})\tilde H_n(x,y)$
%$$=\sum_{i=1}^{n-1}K(\frac{x-X_{ih_n}}{b_{1,n}},
%\frac{y-Y_{ih}}{b_{2n}})\big(\frac{\sigma
%(X_{ih_n})(W_{(i+1)h_n}-W_{ih_n})}{h_n}+g(X_{ih_n},Y_{ih_n})
%\big).$$Thus
\[
\tilde H_n(x,y)=\mathcal I_{1,n}+\mathcal
I_{2n},
\]
where
\[
\mathcal{I}_{1n}=\frac{1}{(n-1)b^d_{1,n}b^d_{2,n}}\sum_{i=1}^{n-1}K
\biggl(\frac
{x-X_{ih_n}}{b_{1,n}},\frac{y-Y_{ih}}{b_{2,n}}\biggr) \biggl(\frac{s
(X_{ih_n})(W_{(i+1)h_n}-W_{ih_n})}{h_n}
\biggr),
\]
and
\[
\mathcal I_{2n}=\frac{1}{(n-1)b^d_{1,n}b^d_{2,n}}\sum_{i=1}^{n-1}K
\biggl(\frac
{x-X_{ih_n}}{b_{1,n}},\frac{y-Y_{ih}}{b_{2,n}}\biggr)g(X_{ih_n},Y_{ih_n}).
\]
We note that
\[
\bigl((n-1)b^d_{1,n}b^d_{2,n}h_n
\bigr)\operatorname{Var}(\mathcal I_{1n})\to s(x)s^*(x)p^\beta_s(x,y)
\int K^2(u,v)\,du\,dv,
\]
and, using the $\alpha$-mixing properties, that
\[
\bigl((n-1)b^d_{1,n}b^d_{2,n} \bigr)
\operatorname{Var}(\mathcal I_{2n})\to g^2(x,y)p^\beta_s(x,y)
\int K^2(u,v)\,du\,dv.
\]
These two last results entail
\[
\lim_{n\to\infty} \bigl((n-1)b^d_{1,n}b^d_{2,n}h_n
\bigr)\operatorname {Var}\bigl(\tilde H_n(x,y)\bigr)\to
s(x)s^2(x)p^\beta_s(x,y)\int
K^2(u,v)\,du\,dv.
\]
Using approximation (\ref{approx}) then allows to derive
\[
\lim_{n\to\infty} \bigl((n-1)b^d_{1,n}b^d_{2,n}h_n
\bigr)\operatorname{Var}\bigl( H_n(x,y)\bigr)\to
\sigma^2(x)p^\beta_s(x,y)\int
K^2(u,v)\,du\,dv.
\]
In this manner, we get
\[
\hat g_n(x,y)=\frac{H_n(x,y)}{\tilde p_s(x,y)}\stackrel{\mathbb P_z}
\to g(x,y).
\]
%
%In particular an estimator of $\beta$ can be readily obtained by
%defining
%\begin{eqnarray*}\label{beta}\hat\beta_n=\frac{\nabla_x
%p_n(x,0)}{H_n(x,0)}\stackrel{\mathbb P_z}\to\beta.\end{eqnarray*}
This result is also true for $d>1$. Thus, using the relation
$g(x,y)+\nabla V(x)=-s(x)s^*(x)y$, we get $\hat g_n(x,y)-\hat
g_n(x,0)\stackrel{\mathbb P_z}\to-s(x)s^*(x)y$, as far as $- \langle\hat
g_n(x,\mathbf e_i)-\hat g_n(x,0),\mathbf e_j \rangle\stackrel{\mathbb P_z}\to
(s(x)s^*(x))_{ij}$.
%\end{eqnarray*}}

%s6 #&#
\section{Examples and numerical simulation results}\label{numexp}
In this section, we want to illustrate some of the main results of the paper.
We start with the It\^o stochastic differential equation defined by
(\ref{WuGirsanov}):
\[
\cases{
dX_t =  Y_t \,dt,
\cr\vspace*{2pt}
dY_t  =  \sigma(X_t,Y_t)
\,dW_t-\bigl(c(X_t,Y_t)Y_t+\nabla
V(X_t)\bigr) \,dt.}
\]
More precisely, we first consider an harmonic oscillator that is driven
by a white noise forcing:
%
%e6.1 #&#
\begin{equation}
\label{OA}
\cases{
dX_t = Y_t
\,dt,
\cr
dY_t  = \sigma \,dW_t-(\kappa Y_t+D
X_t) \,dt, }
\end{equation}
with $\kappa>0$ and $D>0$.
For this model, we know that the stationary distribution is Gaussian,
with mean zero and an explicit variance matrix given in, for example,~\cite{Gardiner}.

For this example, the diffusion term is constant, equal to $\sigma$.
Recall that the infill estimator with $T=1$ is defined by \eqref{estim}:
\[
\hat{\sigma}_n^2=\frac{1}{[{1}/({2h_n})]-1} \frac
{3}{2h_n^3}
\sum_{p=1}^{[{1}/({2h_n})]-1} (X_{(2p+1)h_n}-2X_{2ph_n}+X_{(2p-1)h_n}
)^2.
\]
As the model satisfies assumptions $\mathcal{H}_0$, $\mathcal{H}_1$ and
$\mathcal{H}_2$, we know from Corollary~\ref{theoconstantcVnuls} that
if $h_n \displaystyle{\mathop{\longrightarrow}_{n \rightarrow+ \infty}} 0$, starting from any
initial point $z=(x,y)$,
\[
\sqrt{\frac{1}{2h_n}} \bigl(\hat{\sigma}_n^2-
\sigma^2 \bigr) \mathop{\longrightarrow}_{n \rightarrow+ \infty}^{\mathcal{S}}
\mathcal{N} \bigl(0,2 \sigma^4\bigr).
\]
A $95 \%$ asymptotic confidence interval for $\sigma^2$ is thus
defined as
\[
\mathbf{CI}_{95 \%}\bigl(\sigma^2\bigr)= \bigl[ \hat{
\sigma}_n^2 - 1.96 \sqrt{2} \hat{\sigma}_n^2
\sqrt{2h_n}, \hat{\sigma}_n^2 + 1.96 \sqrt
{2} \hat{\sigma}_n^2 \sqrt{2h_n} \bigr].
\]

In the following, we approximate the solution of \eqref{OA} by an
explicit Euler scheme. We choose $h_n= n^{-\gamma}$, $\gamma>0$,
$\kappa= 2$ and $D=2$.
Then, for different values of $n$ and $\gamma$, we compute $M=1000$
realizations of $\hat{\sigma}_n^2$.
On these $M$ realizations, we compute the empirical relative mean
squared error defined by $\mathrm{RMSE}=\frac{1}M   \sum_{j=1}^M
(({\hat{\sigma}_n^{2,j}-\sigma^2})/{\sigma^2} )^2$, as far as the empirical
coverage of the $95 \%$ confidence interval defined as $\mathrm{ECOV}=\frac{1}M
\sum_{j=1}^M \mathbf{1}_{\sigma^2 \in\mathbf{CI}^j_{95 \%}(\sigma^2)}$.
The results are summarized in Table~\ref{table1} below.

%t1 #&#
\begin{table}[b]
\tablewidth=240pt
\caption{Infill estimation, empirical relative mean squared error
(RMSE) and empirical coverage (ECOV) of the $95$\% confidence interval
with $h_n= n^{-\gamma}$, $M=1000$ realizations of the estimator, and
for different values of $n$, $\gamma$ and $\sigma$}
\label{table1}
\begin{tabular*}{240pt}{@{\extracolsep{\fill}}lcccc@{}}
\hline
$\bolds{\sigma}$ & $\bolds{\gamma}$ & $\bolds{n}$ & \textbf{RMSE} & \textbf{ECOV} \\
\hline
$1$ & $0.5$ & $100$ & $0.47$\phantom{0} &$0.85$ \\
$1$ & $0.5$ & $1000$ & $0.13$\phantom{0}&$0.92$ \\
$1$ & $0.5$ & $10^4$ & $0.04$\phantom{0}&$0.93$ \\
$1$ & $0.7$ & $100$ & $0.19$\phantom{0}&$0.90$ \\
$1$ & $0.7$ & $1000$ & $0.03$\phantom{0}&$0.94$ \\
$1$ & $0.7$ & $10^4$ & $0.006$&$0.95$ \\
$2$ & $0.5$ & $100$ & $2.03$\phantom{0}&$0.86$ \\
$2$ & $0.5$ & $1000$ & $0.53$\phantom{0}&$0.91$ \\
$2$ & $0.5$ & $10^4$ & $0.15$\phantom{0}&$0.94$ \\
$2$ & $0.7$ & $100$ & $0.72$\phantom{0}&$0.91$ \\
$2$ & $0.7$ & $1000$ & $0.13$\phantom{0}&$0.94$ \\
$2$ & $0.7$ & $10^4$ & $0.02$\phantom{0}&$0.95$ \\
\hline
\end{tabular*}
\end{table}

As expected, the more $\gamma$ is high, the more fast is the
convergence. The speed of convergence also depends (through a constant
term in the asymptotic variance) on the unknown value of $\sigma^2$.

We now consider for the same model the infinite-horizon
estimation.
Model~\eqref{OA} satisfies assumptions $\mathcal{H}_0$ up to
$\mathcal{H}_4$. Thus, if $h_n
 \displaystyle{\mathop{\longrightarrow}_{n \rightarrow+ \infty}} 0$,\break  $nh_n
\displaystyle{\mathop{\longrightarrow}_{n \rightarrow+ \infty}} {+}\infty$ 
and $nh_n^2
\displaystyle{\mathop{\longrightarrow}_{n \rightarrow+ \infty}} 0$, then through Theorem~\ref
{cvergconstant}, we have
\[
\sqrt{n} \bigl( \mathcal{K}_n-\sigma^2 \bigr) \mathop{
\longrightarrow}_{n \rightarrow+ \infty}^{\mathcal{D}} \mathcal{N}\bigl(0, 2
\sigma^4\bigr),
\]
with $\mathcal{K}_n= \frac{3}{2}     \frac{1}{(n-1)h_n^3}
\sum_{p=1}^{n-1} (X_{(2p+1)h_n}-2X_{2ph_n}+X_{(2p-1)h_n} )^2   $.
A $95 \%$ asymptotic confidence interval for $\sigma^2$ is thus
defined as
\[
\mathbf{CI}_{95 \%}\bigl(\sigma^2\bigr)= \biggl[
\mathcal{K}_n - 1.96 \frac
{\sqrt
{2} \mathcal{K}_n}{\sqrt{n}}, \mathcal{K}_n +
1.96 \frac
{\sqrt
{2} \mathcal{K}_n}{\sqrt{n}} \biggr].
\]

In the following, we approximate the solution of \eqref{OA} by an
explicit Euler scheme. We choose $h_n= n^{-\gamma}$, $\gamma>0$,
$\kappa= 2$ and $D=2$.
Then, for different values of $n$ and $\gamma$, we compute\vspace*{1pt} $M=1000$
realizations of $\mathcal{K}_n$.
On these $M$ realizations, we compute the empirical relative mean
squared error defined by $\mathrm{RMSE}=\frac{1}M   \sum_{j=1}^M  (\frac
{\mathcal{K}_n^{j}-\sigma^2}{\sigma^2} )^2$, as far as the
empirical coverage of the $95 \%$ confidence interval defined as
$\mathrm{ECOV}=\frac{1}M \sum_{j=1}^M \mathbf{1}_{\sigma^2 \in\mathbf
{CI}^j_{95 \%
}(\sigma^2)}$. The results are summarized in Table~\ref{table2} below.

%t2 #&#
\begin{table}[b]
\tablewidth=240pt
\caption{Infinite-horizon estimation, empirical relative mean squared
error (RMSE) and empirical coverage (ECOV) of the $95\%$ confidence
interval with $h_n= n^{-\gamma}$, $M=1000$ realizations of the
estimator, and for different values of $n$, $\gamma$ and $\sigma$}
\label{table2}
\begin{tabular*}{240pt}{@{\extracolsep{\fill}}lcccc@{}}\hline
$\bolds{\sigma}$ & $\bolds{\gamma}$ & $\bolds{n}$ & \textbf{RMSE} & \textbf{ECOV} \\
\hline
$1$ & $0.5$ & \phantom{0}$100$ & $0.022$&$0.890$ \\
$1$ & $0.5$ & \phantom{0}$500$ & $0.005$&$0.917$ \\
$1$ & $0.5$ & $1000$ & $0.002$&$0.923$ \\
$1$ & $0.7$ & \phantom{0}$100$ & $0.019$&$0.942$ \\
$1$ & $0.7$ &\phantom{0}$500$ & $0.004$&$0.947$ \\
$1$ & $0.7$ & $1000$ & $0.002$&$0.949$ \\
$2$ & $0.5$ & \phantom{0}$100$ & $0.084$&$0.892$ \\
$2$ & $0.5$ & \phantom{0}$500$ & $0.017$&$0.921$ \\
$2$ & $0.5$ & $1000$ & $0.008$&$0.933$ \\
$2$ & $0.7$ & \phantom{0}$100$ & $0.085$&$0.926$ \\
$2$ & $0.7$ & \phantom{0}$500$ & $0.018$&$0.936$ \\
$2$ & $0.7$ & $1000$ & $0.008$&$0.947$ \\
\hline
\end{tabular*}
\end{table}

As expected, we observe that the rate of convergence does not depend on
$\gamma$. The result of Theorem~\ref{cvergconstant} has to be compared
to the one in Theorem~2 in \cite{samthi}. In \cite{samthi}, the
estimator is obtained by minimizing a contrast.
More precisely, the authors in \cite{samthi} define the contrast to
minimize as
\[
\mathcal{L}_n\bigl(\sigma^2\bigr) =\sum
_{p=1}^{n-2} \frac{3}2 \frac{
(X_{(p+1)h_n}-2X_{ph_n}+X_{(p-1)h_n} )^2}{h_n^3 \sigma
^2}+(n-2)
\log \bigl(\sigma^2\bigr),
\]
and they obtain
\[
\tilde{\sigma}_n^2=\frac{3}{2} \frac{1}{n-2}
\sum_{p=1}^{n-2} \frac
{ (X_{(p+1)h_n}-2X_{ph_n}+X_{(p-1)h_n} )^2}{h_n^3}.
\]
They\vspace*{1pt} obtain the same rate of convergence but with the asymptotic
variance equal to $\frac{9}4   \sigma^4$. Our definition \eqref
{acdouble} of the double increment of $X$, which is different from
theirs, allows to recover the asymptotic variance $2   \sigma^4$ they
get for the case of complete observations. In the present paper, we do
not study the optimality of the estimators. It is naturally a very
interesting problem, which, for the model under study is still open.

We now consider a variant of model \eqref{OA} in which we consider a
diffusion term which is nonconstant. It may indeed be interesting in
the applications to choose a position-dependent diffusion term, for
example, to restrict the action of a thermostat to the boundaries only.

More precisely, we consider the following model:
%
%e6.2 #&#
\begin{equation}
\label{OApert}
\qquad\cases{
 dX_t = Y_t
\,dt,
\vspace*{5pt}\cr
\displaystyle dY_t  = \bigl(2 \beta^{-1} \bigr)^{1/2} \exp
\biggl(\frac
{-1}{X_t^2+1} \biggr) \,dW_t- \biggl(\exp \biggl(
\frac
{-2}{X_t^2+1} \biggr) Y_t+ \sin(X_t) \biggr) \,dt.
}\hspace*{-6pt}
\end{equation}

Model \eqref{OApert} is of the form of (\ref{LangevinEquations}).
It satisfies the Einstein's fluctuation-dissipation relation discussed
in Section~\ref{Langevin}.
The potential\vadjust{\goodbreak} is the periodic potential $V(x)=-\cos(x)$ and the
diffusion term is mainly active at the boundaries $s^2(x)=\exp
(\frac{-2}{x^2+1} )$, satisfying\vspace*{1pt} however assumption $\mathcal{H}_0$.

The invariant density is known for that model, but it is possible to
apply the Kernel estimation procedure proposed in \cite{CLP} to
estimate it. In Figures~\ref{densuniv}, \ref{densbiv} below, we chose
$\beta=2$, the Epanechnikov kernel, the bandwidths
$b_{1,n}=b_{2,n}=n^{-0.2}$, and the discretization step $h_n=n^{-0.30}$
with $n=10^5$.

%f1 #&#
\begin{figure}[t]

\includegraphics{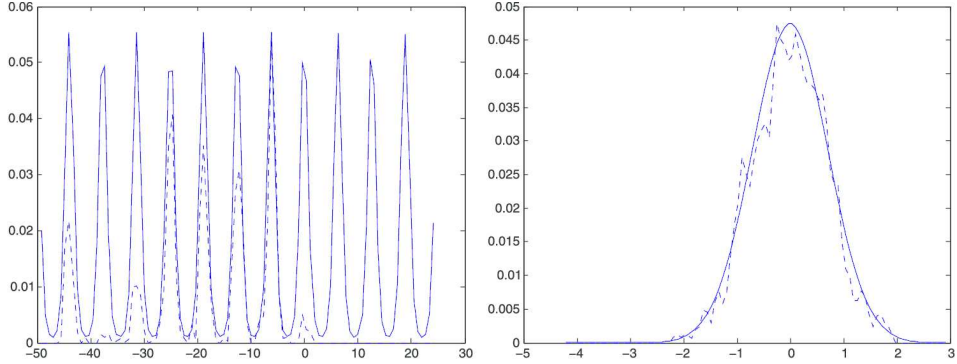}

\caption{Estimated (dashed line) and theoretical (solid line)
invariant density for the position (left) and for the velocity (right),
$\beta=2$, $n=10^5$.}\vspace*{-6pt}
\label{densuniv}
\end{figure}
%f2 #&#
\begin{figure}

\includegraphics{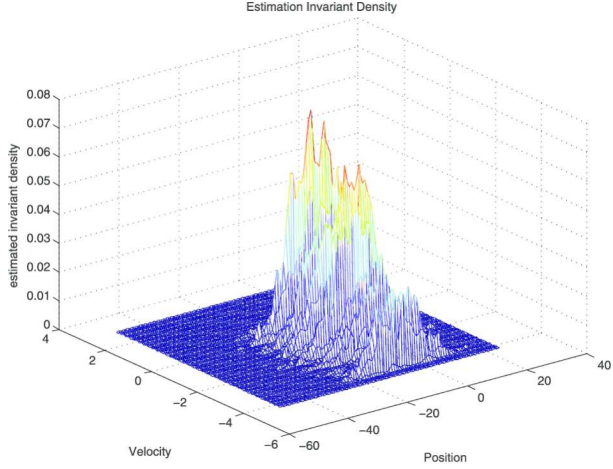}

\caption{Estimated bivariate invariant density, $n=10^6$, $\beta=2$,
$n=10^5$.}
\label{densbiv}
\end{figure}

%\vspace{-4em}
%
%f3 #&#
\begin{figure}[b]

\includegraphics{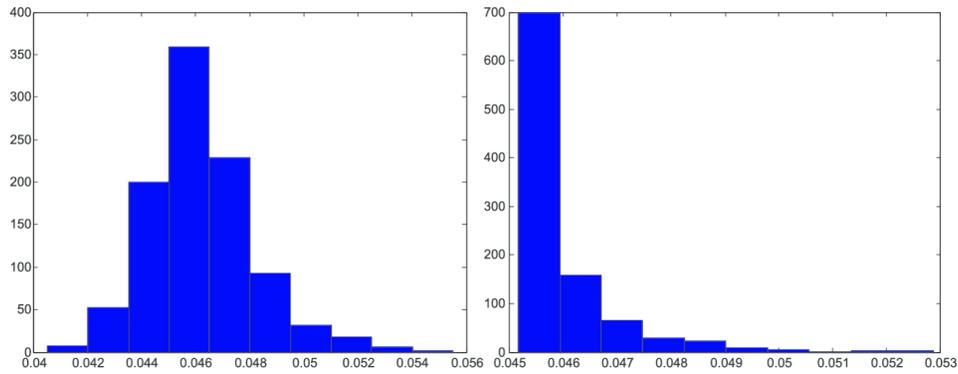}

\caption{Histograms on $M=1000$ realizations of the estimator (left)
and of the limit integral (right), $n=10^5$, $\beta=2$, $h_n=n^{-0.7}$.}
\label{fig1}
\end{figure}

We are only considering then the infill estimation. In that case, the
infill estimator is defined as
\begin{eqnarray*}
&&\mathcal{QV}_{h_n}(1)=\frac{1}{h_n^2} \sum
_{p=1}^{ [{1}/({2h_n}) ]-1} (X_{2p+1)h_n}-2X_{2ph_n}+X_{(2p-1)h_n}
)^2.
\end{eqnarray*}

Thus, if $h_n \displaystyle{\mathop{\longrightarrow}_{n \rightarrow+ \infty}} 0$, we get from
Theorem~\ref{thmmainvarinfill}
\[
\sqrt{\frac{1}{h_n}} \biggl(\mathcal{QV}_{h_n}(1)-
\frac
{2}{3\beta} \int_0^1 \exp \biggl(
\frac{-2}{X_s^2+1} \biggr)\,ds \biggr)\mathop{\longrightarrow}_{n \rightarrow+ \infty}^{\mathcal S}
\frac{4}{3\beta} \int_0^1 \exp \biggl(
\frac{-2}{X_s^2+1} \biggr) \,d\tilde W_s,
\]
where $( \tilde{W}_t,   t \in[0,T] )$ is a Wiener
process independent of the initial Wiener process $W_{.}$, with
variance equal to 2.

In the following, we choose $h_n= n^{-\gamma}$ with $\gamma=0.7$. We
compute $M=1000$ realizations of the estimator $\mathcal{QV}_{h_n}(1)$
and $M=1000$ realizations of the limit $\frac{2}{3\beta}   \int_0^1
\exp (\frac{-2}{X_s^2+1} )\,ds$. This integral is approximated
by a quadrature formula with the rectangle rule.

We consider the case $n=10^5$, $\beta=2$. We compute the empirical
relative mean squared error (RMSE) and we draw (see Figure~\ref{fig1})
both the histogram of the estimator and the one of the limit integral
for the $M=1000$ realizations.

We get for both cases $\mathrm{RMSE}= 0.0024$.

%
%\begin{figure}[ht!]\footnotesize
%\centering
%\subfigure[]{\includegraphics[scale=0.27]{est2.pdf}}
%~~~~~~
%\subfigure[]{\includegraphics[scale=0.27]{approx2.pdf}}
%\caption{Histograms on $M=1000$ realizations of the estimator and of
%the limit integral, $n=10^4$, $\sigma=2$, $h_n=n^{-0.7}$.}
%\label{fig2}
%\end{figure}

The histograms in Figures~\ref{fig1} (left and right) have
similarities. However, we note that we have an important boundary
effect for the lower tail in both cases, probably due to the
approximation of the limit integral by a quadrature rule.

\section*{Acknowledgements}
We want to heartily acknowledge an anonymous referee for  comments
on a first version of the paper, especially for suggesting to describe
how our techniques can apply to the Langevin dynamics in Section~\ref{Langevin}.

% imsref loaded by daiva.urboniene, 2015-08-06 08:42:08

%\begin{appendix}
%\section{}
%\end{appendix}

%\section*{Acknowledgments}

%\begin{supplement}[id=suppA]
%\slink[doi]{10.1214/00-AAPXXXXSUPP} %[doi,text={...}] - jei reikia
%suskaldyti doi
%\sfilename{aapXXXX\_supp.pdf}
%\end{supplement}

%\begin{thebibliography}{99}
%\bibitem[\protect\citeauthoryear{}{}]{r1}
%\bibitem{r1}
%\end{thebibliography}

\printaddresses
\end{document}